\documentclass[a4paper]{article}

\usepackage{amsmath,amsthm,amssymb}
\usepackage[matrix,curve,arrow,tips,frame]{xy}
\usepackage{tabularx}
\usepackage{enumerate,epic,eepic}

   \newtheorem{thm}{Theorem}[section]
   \newtheorem{lem}[thm]{Lemma}
   \newtheorem{cor}[thm]{Corollary}
   \newtheorem{prop}[thm]{Proposition}
   \newtheorem{defn}{Definition}[section]
   \newtheorem{exmp}{Example}[section]
   \newtheorem{rem}{Remark}[section]
   \newtheorem{pf}{Proof}

\allowdisplaybreaks[1]

\numberwithin{equation}{section}

\DeclareFontEncoding{OT2}{}{}
\DeclareFontSubstitution{OT2}{wncyr}{m}{n}
\DeclareSymbolFont{cyss}{OT2}{wncyss}{m}{n}
\DeclareMathSymbol{\sh}{\mathbin}{cyss}{`x}

\title{Finite automata and relations of multiple zeta values}
\author{Sinya KITANI,\quad Eiki SAWADA, \quad Kimio UENO}
\date{\today}

\begin{document}
\maketitle

\begin{abstract}
The theory of finite automata applies to the study on relations of multiple zeta
values.
\end{abstract}

\section{Introduction}

\subsection{Multiple zeta values and Zagier - Broadhurst's formula}

Multiple zeta values (MZVs, for short) are real numbers defined by
\begin{align}\label{MZV}
 \zeta(k_1,k_2,\ldots,k_n) = \sum_{m_1>m_2>\cdots>m_n>0} \ 
         \frac{1}{m_1^{k_1}m_2^{k_2} \cdots m_n^{k_n}}
\end{align}
where $k_1,k_2,\ldots,k_n \in \mathbb{Z}_{\geq 1}$, and $k_1 \geq 2$.
We call $k=k_1+\cdots+k_n$ the weight and n the depth of \eqref{MZV}.

The following formula was conjectured by Zagier \cite{Z}, and proved by Broadhurst \cite{BB},
\cite{BBBL}:
\begin{align}\label{Zagier-Broadhurst}
\zeta(\{3,1\}_n) 
    = \zeta(\underbrace{3,1,\ldots,3,1}_{n\,\text{times}})
    =\frac{2 \pi^{4n}}{(4n+2)!}.
\end{align}
The proof refined by Zagier is as follows (cf. \cite{AK}): The equality
\begin{align*}
    \sum_{n=0}^{\infty} Li_{\underbrace{\scriptstyle3,1,\ldots,3,1}_{n\,\text{times}}}(x)\,t^{4n}
     =F\Big( \frac{t}{1+i},\frac{-t}{1+i};1;x \Big)\,F\Big( \frac{t}{1-i},\frac{-t}{1-i};1;x \Big)
\end{align*}
can be proved by showing that the both sides is annihilated by the differential operator
\begin{align*}
      \Big( (1-x)\frac{d}{dx} \Big)^2 \Big( x\frac{d}{dx} \Big)^2\,-\,t^4.
\end{align*}
Here, for $k_1,\ldots,k_n \in \mathbb{Z}_{\geq1}$,
\begin{align}
       Li_{k_1,\ldots,k_n}(z) = \sum_{m_1>\cdots>m_n>0}
                          \frac{z^{m_1}}{m_1^{k_1} \cdots m_n^{k_n}}
\end{align}
are multiple polylogarithms of one variable (MPLs, for short). Hence we have
\begin{align*}
  \sum_{n=0}^{\infty} \zeta(\underbrace{3,1,\ldots,3,1}_{n\,\text{times}})\,t^{4n}
      =F\Big( \frac{t}{1+i},\frac{-t}{1+i};1;1 \Big)\,
                       F\Big( \frac{t}{1-i},\frac{-t}{1-i};1;1 \Big).
\end{align*}
Applying the formula
\begin{align*}
  F(a,-a;1;1)=\frac{1}{\Gamma(1-a)\Gamma(1+a)}=\frac{\sin \pi a}{\pi a},
\end{align*}
we have
\begin{align*}
  F\Big( \frac{t}{1+i},\frac{-t}{1+i};1;1 \Big)\, F\Big( \frac{t}{1-i},\frac{-t}{1-i};1;1 \Big)
  =
  \sum_{n=0}^{\infty}\frac{2\pi^{4n}}{(4n+2)!}\,t^{4n}.
\end{align*}

\subsection{Waldschmidt's idea}

The original proof by \cite{BB}, \cite{BBBL} is based upon some combinatorics, 
while Waldschmidt gave more simple one by means of the idea ``finite automata''. 

Let $\mathfrak{H}=\mathbb{Q}\langle x,y \rangle$ be a $\mathbb{Q}$-algebra 
of polynomials of non-commutaitive variables $x$ and $y$, and a subalgebra
$\mathfrak{H}^0 = \mathbb{Q}\cdot1 + x\mathfrak{H}y$. Set $z_k=x^{k-1}y \ (k \geq 1)$. We define
a $\mathbb{Q}$-linear map $Z \,:\, \mathfrak{H}^0 \longrightarrow \mathbb{R}$ \ by
\begin{align}
     Z(z_1 z_2 \cdots z_n) = \zeta(k_1,k_2,\ldots,k_n) \quad
                  ( k_1,k_2,\ldots,k_n \in \mathbb{Z}_{\geq 1}  \ \text{and} \ k_1 \geq 2).
\end{align}
which coressponds to the iterative integral representation;

\begin{align}
 \zeta(k_1,k_2,\ldots,k_n) = 
    \int_0^1 \underbrace{\frac{dt}{t} \cdots \frac{dt}{t}}_{k_1-1} \frac{dt}{1-t}
     \cdots
     \underbrace{\frac{dt}{t} \cdots \frac{dt}{t}}_{k_n-1} \frac{dt}{1-t}.
\end{align}

The shuffle product $\sh$ is a notion abstracting the product of iterated integrals.

\begin{defn}
$\sh : \mathfrak{H} \times \mathfrak{H} \longrightarrow \mathfrak{H}$ is a 
$\mathbb{Q}$-bilinear operation satisfying
the following conditions:
\begin{enumerate}[$(i)$]
\item For $w \in \mathfrak{H}$, $w \sh 1 = 1 \sh w = w$.
\item Let $u_1, u_2$ be $x$ or $y$. Then, for any words $w_1, w_2 \in \mathfrak{H}$,
\begin{align}
           (u_1w_1) \sh (u_2w_2) = u_1(w_1 \sh u_2w_2) + u_2(u_1w_1 \sh w_2).
\end{align}
\end{enumerate}
\end{defn} 

It is known that $(\mathfrak{H},\sh)$ is a commutative algebra (cf. \cite{AK}), and 
$(\mathfrak{H}^0,\sh)$ is a commutative subalgebra, furthermore, $Z$ is an algebra homomorphism,
that is,
\begin{align}\label{eq:homsh}
    Z(w_1 \sh w_2) = Z(w_1)Z(w_2),   \quad  (w_1,w_2 \in \mathfrak{H}^0).
\end{align}

Let $\widehat{\mathfrak{H}}= \mathbb{Q}\langle\langle x,y \rangle\rangle$ be a 
$\mathbb{Q}$-algebra of formal power series of $x, \ y$ which can be viewed as the dual of 
$\mathfrak{H}$. For $w \in x\mathfrak{H}+y\mathfrak{H}$, set
\begin{align}\label{closure}
           w^{\star}=\sum_{n=0}^{\infty}w^n
\end{align} 
which is called Kleene's closure of $w$. This is the inverse element of $w$ in 
$\widehat{\mathfrak{H}}$:
\begin{align}
    w^{\star}(1-w) \,=\, (1-w)w^{\star} \,=\, 1.
\end{align}
Using the idea of finite automata, Waldschmidt \cite{W} showed that
\begin{align}\label{Waldschmidt}
   (xy)^{\star} \sh (-xy)^{\star} = (-4x^2y^2)^{\star}.
\end{align}
We should note that the shuffle product preserves the weight of words.
Transforming the both sides of \eqref{Waldschmidt} via $Z$ to MZVs, one obtaines 
\begin{align*}
     (-4)^n\zeta(\{3,1\}_n) = \sum_{p=0}^{2n}(-1)^{2n-p}\zeta(\{2\}_p)\zeta(\{2\}_{2n-p}).
\end{align*}
The RHS above equals to $(-1)^n \zeta(\{4\}_n)$  \ (see Corolally \ref{cor:zhkk}
\ \eqref{eq:zeta2kn}), and
\begin{align}\label{val:zeta4}
     \zeta(\{4\}_n) = \frac{2\cdot4^n \pi^{4n}}{(4n+2)!},
\end{align}
so that \eqref{Zagier-Broadhurst} is established.

The original proof of \eqref{Waldschmidt} will be reviewed in Section 3 after the preliminaries
on automata theory.

\subsection{The purpose of the paper}

In this paper, we aim at extending the idea of Waldschmidt and
deriving various relations of MZVs.
Our central idea is to associate an "adjacency matrix" 
to each finite automaton. Through adjacency matrices, one can compute the shuffle product
in a combinatorial way. %Our method effectively works even in the case of the shuffle 
%product and the harmonic product of multiple $L$ values, however these subjects will 
%be treated in a forthcoming paper.

This paper is organized as follows: In Section 2, we give preliminaries on the basic concept
of finite automata, in particular, shuffle automata and adjacency matrices. In Section 3,
we review in details the proof of \eqref{Waldschmidt} due to \cite{W}, and give the definition
of the harmonic product of MZVs and show basic formulas which will be used later. 
In Section 4, first we give a direct application of the Waldschmidt formula
\eqref{Waldschmidt}. Next we derive several relations of MZVs by means of 
adjacency matrices. In Section 5 we introduce hamonic automata, and prove the basic formula
presented in Section 3. In Appendix we consider the values of $\zeta(\{2k\}_n) \
(1\leq k \leq 7)$.

In \cite{K} and \cite{S}, the first and the second authors of this article deal with
more shuffle automata and harmonic automata with generalization to multiple L values.
These subjects will be treated in the next paper.
\vspace{3mm}

\noindent
{\large \textbf{Acknowledgement}}

\vspace{3mm}

\noindent
This research started with inspired by the talk of Professor Waldschmidt at the conference 
on ``Zeta Functions, Topology, and Quantum Physics'' held at Kinki University, March 2003.
The authors would like to express deep gratitude to Professor Waldschmidt and the organizing
committee of ZTQ.

The third author is partially supported by JPSP Grant-in-Aid No. 15540050.

\section{Finite automata}

\subsection{Defintion of finite automata}

\begin{defn}\label{F/A}
\begin{enumerate}[$(1)$]
\item A finite automaton over $\mathfrak{H}_{\mathbb{C}}= \mathbb{C} \langle x,y \rangle$ is a
      quintuple $(Q,\Sigma,\delta,q_1,F)$ where
\begin{enumerate}[$(i)$]
\item $Q=\{q_1,q_2,\ldots,q_m\}$ is a set of states;
\item $\Sigma$ is a finite subset of $\mathbb{C}x+\mathbb{C}y$
      which is called the alphabets of the automaton;
\item $\delta : Q \times \Sigma \longrightarrow Q$ \ is the transition function;
\item $q_1$ is the initial state;
\item $F \subset Q$ is the set of final states.
\end{enumerate}
\item If $w=u_1 \cdots u_n \in \mathfrak{H}_{\mathbb{C}} \ (u_1,\ldots,u_n \in \Sigma)$ satisfies
\begin{align*}
    \delta(\cdots(\delta(q_1,u_1) \cdots u_n)=q_i \in F,
\end{align*}
then we say that the word $w$ is accepted by the automaton. 
\item If, in the automaton, the states transit like
\begin{align*}
    q_i=q_{i_0} \rightarrow q_{i_1} \rightarrow \cdots q_{i_{n-1}} \rightarrow q_{i_n}=q_j,
\end{align*}
then we call this trajectry a path of length n from $q_i$ to $q_j$. 
If the word $w$ is accepted by a path of length n, then we call $w$ a word of
length n, and denote it by $l(w)=n$.
\item If, for $w \in \widehat{\mathfrak{H}}_{\mathbb{C}} = 
\mathbb{C} \langle\langle x,y \rangle\rangle$, 
there exists an automaton $M=(Q,\Sigma,\delta,q_1,F)$ such that
\begin{align*}
      w = & \sum \,(\text{a word accepted by the automaton}) \hspace{2cm}\\
             & \hspace{1cm}   \times (\text{the number of paths for which the word is accepeted}), 
\end{align*}
then we say that the element $w=w(M)$ is accepted by the automaton $M$.
\end{enumerate}
\end{defn}

\begin{rem}
In the ordinary theory of automata (cf. \cite{An},\cite{HMU}), what is 
accepted by an automaton is a ``language'', but in our theory, it is an element in 
$\widehat{\mathfrak{H}}_{\mathbb{C}}$.
\end{rem}

\begin{exmp}
Let $M=(Q,\Sigma,\delta,q_1,F)$ be
\begin{align*}
     & Q=\{q_1,q_2,q_3\}, \\
     & \Sigma=\{x,y,-y\},\\
     & \delta: \ \delta(q_1,x)=q_2, \ \delta(q_1,-y)=q_3, \ \delta(q_2,x)=q_3, \ 
               \delta(q_3,y)=q_1, \\
     & F=\{q_3\}.
\end{align*}
The transition diagram of the automaton is
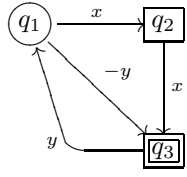
\[\UseTips
\xymatrix @+3mm {
 *++[o][F]{q_1} \ar[r]^{x} \ar[dr]^{-y}
& *+[F]{q_2} \ar[d]^{x}
 \\
 & *+[F=]{q_3} \ar `l_lu[ul]+/d1.1cm/`d_d[ul]_{y} [ul]
}
\]
where the initial state is surrounded by a circle and the final state by a double
rectangular. For this automaton, $w(M)=(x^2-y)(yx^2-y^2)^{\star}$.
\end{exmp}

Now we present addition of elements, multiplication of elements by scalors,
concatenation of elements in terms of automata. 
Let $w_1=u_1u_2 \cdots u_m, \ w_2=v_1v_2\cdots v_n \in \mathfrak{H}
\ (u_1,\ldots,u_m,v_1,\ldots,v_n \in \Sigma)$.
These are represented by the following automata:
\begin{align}
\UseTips
\xymatrix @-1mm {
 w_1: 
& *++[o][F]{q_1} \ar[r]_{u_1}
& *+[F]{q_2} \ar[r]_{u_2}
& \cdots \ar[r]_{u_m}
& *+[F=]{q_{m+1}}
}\ , 
\end{align}
\begin{align}
\UseTips
\xymatrix @-1mm { w_2:
& *++[o][F]{q_1} \ar[r]_{v_1}
& *+[F]{q_2} \ar[r]_{v_2}
& \cdots \ar[r]_{v_n}
& *+[F=]{q_{n+1}}
}\ .
\end{align}
Then the sum $w_1+w_2$ and the concatenation $w_1w_2$ are represented by
\begin{align}\label{addition}
\UseTips
\xymatrix @-1mm {
 w_1+w_2: 
& *++[o][F]{q_1} \ar[r]_{u_1} \ar[rd]_{v_1}
& *+[F]{q_2} \ar[r]_{u_2}
& \cdots \ar[r]_{u_m}
& *+[F=]{q_{m+1}}
\\
& 
& *+[F]{q_{m+2}} \ar[r]_{v_2}
& \cdots \ar[r]_{v_n}
& *+[F=]{q_{m+n+1}}
}\ ,
\end{align}
\begin{align}
\UseTips
\xymatrix @-1mm {
 w_1w_2: 
& *++[o][F]{q_1} \ar[r]_{u_1}
& *+[F]{q_2} \ar[r]_{u_2}
& \cdots \ar[r]_{u_m}
& *+[F]{q_{m+1}} \ar[r]_{v_1}
& *+[F]{q_{m+2}} \ar[r]_{v_2}
& \cdots \ar[r]_{v_n}
& *+[F=]{q_{m+n+1}}
} \ .
\end{align}
The scalor multiplication ($k \in K$) and Kleene's closure of $w_1$ are represented by 
\begin{align}
\UseTips
\xymatrix @-1mm {
 kw_1: 
& *++[o][F]{q_1} \ar[r]_{u_1}
& *+[F]{q_2} \ar[r]_{u_2}
& \cdots \ar[r]_{ku_m}
& *+[F=]{q_{m+1}}
}\ ,
\end{align}

\begin{align}\label{fig:closure}
\UseTips
\xymatrix @-1mm {
w_1^{\star}:
& *++[o][F=]{q_1} \ar[r]_{u_1}
& *+[F]{q_2} \ar[r]_{u_2}
& \cdots \ar[r]_{u_{m-1}}
& *+[F]{q_{m}} \ar `dl_l[lll]+/d6mm/`l_ul[lll]^{u_m} [lll]
}\ .
\end{align}
In \eqref{fig:closure}, ``the initial state $=$ the final state'' is surrounded by a double circle.
We should note that the automata representing the addition and the
scalor multiplication is not unique. For example, if $u_1=v_1=x$, then 
$w_1+w_2$ is also represented by
\begin{align*}
\UseTips
\xymatrix @-1mm {
 w_1+w_2: 
& *++[o][F]{q_1} \ar[r]_{x}
& *+[F]{q_2} \ar[r]_{u_2} \ar[rd]_{v_2}
& *+[F]{q_3} \ar[r]_{u_3}
& \cdots \ar[r]_{u_m}
& *+[F=]{q_{m+1}}
\\
&
& 
& *+[F]{q_{m+2}} \ar[r]_{v_3}
& \cdots \ar[r]_{v_n}
& *+[F=]{q_{m+n}}
}\ .
\end{align*}
In this automaton, the number of the states is less than that of \eqref{addition}.
To reduce the number of the states in an automaton is important for easy computation of 
the words accepted. %This is crucial for using ``adjacency matrice'' which will be defined later.

\subsection{Shuffle automata and adjacency matrices}

Automata reprsenting the shuffle product of words and of Kleene's closure of words
are considered in \cite{HMU} and \cite{W}.

\begin{prop}
The shuffle product $w_1 \sh w_2$ and $(w_1)^{\star} \sh (w_2)^{\star}$ are represented by
the transition diagrams \eqref{fig:sh} and \eqref{fig:shstar}, respectively:
\begin{align}\label{fig:sh}
\UseTips
\xymatrix @+2mm {
 *++[o][F]{q_1} \ar[r]_{u_1} \ar[d]_{v_1}
& *+[F]{q_2} \ar[r]_{u_2} \ar[d]_{v_1}
& \cdots \ar[r]_{u_m}
& *+[F]{q_{m+1}} \ar[d]_{v_1} 
\\
 *+[F]{q_{(m+1)+1}} \ar[r]_{u_1} \ar[d]_{v_2}
& *+[F]{q_{(m+1)+2}} \ar[r]_{u_2} \ar[d]_{v_2}
& \cdots \ar[r]_{u_m}
& *+[F]{q_{2(m+1)}} \ar[d]_{v_2}
\\
 \vdots \ar[d]_{v_{n-1}}
& \vdots \ar[d]_{v_{n-1}}
& \vdots 
& \vdots \ar[d]_{v_{n-1}}
\\
 *+[F]{q_{(n-1)(m+1)+1}} \ar[r]_{u_1} \ar[d]_{v_n}
& *+[F]{q_{(n-1)(m+1)+2}} \ar[r]_{\phantom{abc}u_2} \ar[d]_{v_n}
& \cdots \ar[r]_{u_m}
& *+[F]{q_{n(m+1)}} \ar[d]_{v_n}
\\
 *+[F]{q_{n(m+1)+1}} \ar[r]_{u_1} 
& *+[F]{q_{n(m+1)+2}} \ar[r]_{u_2} 
& \cdots \ar[r]_{u_m}
& *+[F=]{q_{(n+1)(m+1)}} 
}
\end{align}

\begin{align}\label{fig:shstar}
\UseTips
\xymatrix @+2mm {
 *++[o][F=]{q_1} \ar[r]_{u_1} \ar[d]_{v_1}
& *+[F]{q_2} \ar[r]_{u_2} \ar[d]_{v_1}
& \cdots \ar[r]_{u_{m-1}}
& *+[F]{q_{m}} 
\ar `ul^l[lll]+/u6mm/`l^dl[lll]_{u_m} [lll]  
\ar[d]_{v_1} 
\\
 *+[F]{q_{m+1}} \ar[r]_{u_1} \ar[d]_{v_2}
& *+[F]{q_{m+2}} \ar[r]_{u_2} \ar[d]_{v_2}
& \cdots \ar[r]_{u_{m-1}}
& *+[F]{q_{2m}} 
\ar `ul^l[lll]+/u6mm/`l^dl[lll]_{u_m} [lll]  
\ar[d]_{v_2}
\\
 \vdots \ar[d]_{v_{n-1}}
& \vdots \ar[d]_{v_{n-1}}
& \vdots 
& \vdots \ar[d]_{v_{n-1}}
\\
 *+[F]{q_{(n-1)m+1}} 
 \ar[r]_{u_1}
 \ar `ul_u[uuu]+/l1cm/`u_ru[uuu]^{v_n} [uuu] 
& *+[F]{q_{(n-1)m+2}}
 \ar[r]_{u_2}
 \ar `ul_u[uuu]+/l1cm/`u_ru[uuu]^{v_n} [uuu] 
& \cdots \ar[r]_{u_{m-1}}
& *+[F]{q_{nm}} 
 \ar `dl_l[lll]+/d6mm/`l_ul[lll]^{u_m} [lll]
 \ar `ur^u[uuu]+/r1cm/`u^lu[uuu]_{v_n} [uuu] 
} 
\end{align}
\end{prop}

We call \eqref{fig:sh}, \eqref{fig:shstar} the shuffle automata.
To compute a shuffle product is nothing but computing the words accepted by the 
shuffle automaton. For this end, an adjacency matrix is a useful tool.

\begin{defn}
Let $M=(Q,\Sigma,\delta,q_1,F)$ be a finite automaton. Then we define a matrix $A=A(M)$ by
\begin{align}
  A=(a_{ij}) \quad \text{where}  \quad a_{ij}=\sum_{a\in\Sigma\,;\,\delta(q_i,a)=q_j}a.
\end{align}
\end{defn}
We call $A(M)$ the adjacency matrix of $M$. This satisfies the following property.

\begin{thm}\label{prop:am}
Let us denote the $(i,j)$-entry of $A(M)^n$ by $a_{ij}^{(n)}$. Then $a_{ij}^{(n)}$
is the sum of all the words accepted by paths of length n from $q_i$ to $q_j$.
\end{thm}

\begin{pf}
We prove by induction. The case of $n=1$ is trivial. Assume that the case of $n-1$ holds.
Let $Q=(q_1,q_2,\ldots,q_m)$ be the set of states of $M$. Since $A^n=A^{n-1}A$, we have
\begin{align*}
      a_{ij}^{(n)} = \sum_{k=1}^m a_{ik}^{(n-1)}a_{kj}^{(1)}.
\end{align*}
On the other hand, there exist m strings of paths of length n like
\begin{align*}
    \underbrace{q_i \rightarrow \cdots \rightarrow q_k}_{length=n-1}
     \rightarrow q_j \quad (k=1,2,\ldots,m).
\end{align*}
Hence $a_{ij}^{(n)}$ equals to the sum of all the words accepted by paths of 
length n from $q_i$ to $q_j$.
\hfill$\square$
\end{pf}

\begin{thm}\label{thm:am}
The element $w(M)$ accepted by $M$ is given by
\begin{align}
          w(M) = \sum_{j \,;\, q_j \in F} \left( \sum_{n=0}^{\infty} a_{1j}^{(n)} \right).
\end{align}
(Here, for convenience, we set the word of length 0 to be $1 \in \widehat{\mathfrak{H}}_{\mathbb{C}}$
which is accepted only by the paths from $q_i$ to $q_i$. ) In particular, 
in the case that $F=\{q_1\}$, it is given by
\begin{align}
          w(M) = \sum_{n=0}^{\infty} a_{11}^{(n)}.
\end{align}
\end{thm}

\section{The Waldschmidt formula and harmonic product of MZVs}

\subsection{The original proof of the Waldschmidt formula}
In \cite{W}, Waldschmidt proved the formula \eqref{Waldschmidt} in the following manner:
From \eqref{fig:shstar}, the shuffle automaton of $(xy)^{\star}\sh(-xy)^{\star}$ is represented
by the following transition diagram:
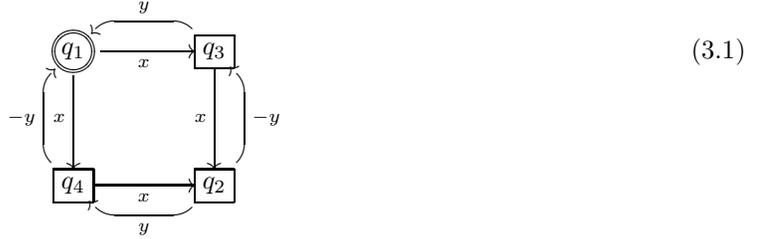
\begin{align}\label{fig:w22} 
\UseTips
\xymatrix @+4mm {
 *++[o][F=]{q_1} \ar[r]_{x} \ar[d]_{x}
& *+[F]{q_3} 
\ar `ul^l[l]+/u4mm/`l^dl[l]_{y} [l]  
\ar[d]_{x} 
\\
 *+[F]{q_4} 
 \ar[r]_{x}
 \ar `ul_u[u]+/l4mm/`u_ru[u]^{-y} [u] 
& *+[F]{q_2} 
 \ar `dl_l[l]+/d4mm/`l_ul[l]^{y} [l]
 \ar `ur^u[u]+/r4mm/`u^lu[u]_{-y} [u] 
 }
\end{align}
Denote this automaton by $M=(Q,\Sigma,\delta,q_1,F)$. Let 
$M_k=(Q,\Sigma,\delta,q_k,F)$ and $S_k$ the element accepted by $M_k$. 
Then we have the following linear recursive equations;
\begin{align}\label{lineareq.}
\left\{
\begin{array}{llll}
 S_1=1+xS_3+xS_4, & \quad S_2=-yS_3+yS_4, \\
 S_3=yS_1+xS_2,  & \quad S_4=-yS_1+xS_2.
\end{array}
\right.
\end{align}
Solving this, we have 
\begin{align}
      S_1 = 1 - 4x^2y^2S_1.
\end{align}
Hence $S_1=(-4x^2y^2)^{\star}$.

For any finite automaton, one has linear recursive equations like \eqref{lineareq.}.
But they cannot be solved (or, it is too difficult to solve) in general. By introducing
adjacency matices, one can avoid such difficulty.

\subsection{Proof of the Waldschmidt formula via an adjacency matrix}
The adjacency matrix $A$ of the automaton \eqref{fig:w22} is
\begin{align}
  A= 
   \begin{bmatrix}
     0 & 0 & x &  x \\
     0 & 0 & y & -y \\
     y & x & 0 &  0 \\
     -y& x & 0 &  0 
   \end{bmatrix}
   .
\end{align}
Let us compute $a_{11}^{(n)}$. We have
\begin{align*}
  A^2 = 
   \begin{bmatrix}
       0     & 2x^2 & 0 & 0 \\
       -2y^2 &   0  & 0 & 0 \\
       0     &   0  & * & * \\
       0     &   0  & * & * 
   \end{bmatrix}
   ,
\end{align*}
where $*$ stand for certain elements of length 2, and 
\begin{align*}
  A^4 = 
   \begin{bmatrix}
    -4x^2y^2 &        0 & 0 & 0 \\
       0     & -4y^2x^2 & 0 & 0 \\
       0     &        0 & * & * \\
       0     &        0 & * & * 
   \end{bmatrix}
   ,
\end{align*}
where $*$ stand for certain elements of length 4. From this computation, we see that
$a_{11}^{(n)} = 0$ unless $n=4k$, and $a_{11}^{(4k)}=(-4x^2y^2)^k$. Thus we obtain
the Waldschmidt formula \eqref{Waldschmidt}.

\subsection{Harmonic product}

Let $k_1, k_2 \geq 2$. Then one can compute product of zeta values 
$\zeta(k_1), \ \zeta(k_2)$ like
\begin{align*}
   \zeta(k_1)\zeta(k_2) & = \left\{ \sum_{m_1>m_2>0}+
               \sum_{m_2>m_1>0}+\sum_{m_1=m_2>0} \right\}\, \frac{1}{m_1^{k_1}m_2^{k_2}} \\
                        & = \zeta(k_1,k_2)+\zeta(k_2,k_1)+\zeta(k_1+k_2).
\end{align*}
Generalizing this, one can introduce the harmonic product $*$ on 
$\mathfrak{H}^1=\mathbb{Q}\cdot1+\mathfrak{H}y$. 
\begin{defn}
$* : \mathfrak{H}^1 \times \mathfrak{H}^1 \longrightarrow \mathfrak{H}^1$
is a $\mathbb{Q}$-bilinear operation satisfying the following conditions:
\begin{enumerate}[$(i)$]
\item For any $w \in \mathfrak{H}^1$, $w * 1 = 1 * w = w$.
\item Let $z_k=x^{k-1}y$ \ $(k=1,2,\ldots)$. For any words $w_1, \ w_2$ in $\mathfrak{H}^1$,
\begin{align}\label{def:harmonic}
       (z_iw_1)*(z_jw_2) = z_i(w_1*(z_jw_2))+z_j((z_iw_1)*w_2)+z_{i+j}(w_1*w_2).
\end{align}
\end{enumerate}
\end{defn}

Hoffman \cite{H} showed that $(\mathfrak{H}^1,*)$ is a commutative algebra generated by 
$z_k=x^{k-1}y \ (k=1,2,\ldots)$, and that
the map $Z : \mathfrak{H}^0 \longrightarrow \mathbb{R}$ is an algebra homomorphism, that is,
\begin{align}\label{eq:homhar}
          Z(w_1*w_2) =Z(w_1)Z(w_2),  \quad  (w_1, \ w_2 \in \mathfrak{H}^0).
\end{align}
By \eqref{eq:homsh} and \eqref{eq:homhar}, we have the finite double shuffle relation:

\begin{prop}[Zagier, cf. \cite{AK}]\label{prop:fdsr}
\begin{align}\label{eq:fdsr}
         Z(w_1*w_2 - w_1 \sh w_2) = 0,  \quad  (w_1, \ w_2 \in \mathfrak{H}^0).
\end{align}
\end{prop}

The harmonic product of the Kleene closure of $z_k$ is computated as follows:
\begin{thm}\label{thm:hkk}
We have
\begin{align}\label{eq:closurehar2}
   (z_k)^{\,\star}*(-z_k)^{\star} = (-z_{2k})^{\star}.
\end{align}
In general, letting $\omega$ be a primitive m-th root of unity,
\begin{align}\label{eq:closureharn}
       (z_k)^{\star}*(\omega z_k)^{\star}*\cdots*(\omega^{m-1}z_k)^{\star} = ((-1)^{m-1}z_{mk})^{\star}.
\end{align}
\end{thm}
%\begin{pf}
%Let $w^{\,+} = \sum_{n=1}^{\infty}w^{\,n} = w\, w^{\,\star}$. Then we have
%\begin{align*}
%     (z_k)^{\star}*(-z_k)^{\star} &= (1+(z_k)^+)*(1+(-z_k)^+) \\
%                                  &= 1 +(z_k)^+ + (-z_k)^+ + (z_k)^+ * (-z_k)^+.
%\end{align*}
%By \eqref{def:harmonic}, we have
%\begin{align*}
%    (z_k)^+ * (-z_k)^+ &= z_k((z_k)^{\star} * (-z_k)^+)-z_k((z_k)^+ * (-z_k)^{\star})
%                                -z_{2k}((z_k)^{\star} * (-z_k)^{\star}) \\
%                       &= z_k(-z_k)^+ - z_k(z_k)^+ -z_{2k}((z_k)^{\star} * (-z_k)^{\star})\\
%                       &= -(-z_k)^+ - (z_k)^+ -z_{2k}((z_k)^{\star} * (-z_k)^{\star}).
%\end{align*}
%Thus we have
%\begin{align*}
%          (z_k)^{\star}*(-z_k)^{\star} = 1- z_{2k}((z_k)^{\star} * (-z_k)^{\star}).
%\end{align*}
%\hfill$\square$
%\end{pf}
We will prove this theorem in Section 5 after introducing the notion of ``harmonic automata''.
Transforming \eqref{eq:closurehar2}, \eqref{eq:closureharn} 
to MZVs via $Z$ yields the following formula:

\begin{cor}\label{cor:zhkk}
We have
\begin{align}\label{eq:zeta2kn}
  \sum_{p=0}^{2n} (-1)^p\zeta(\{ k \}_p) \, \zeta(\{ k \}_{2n-p}) = (-1)^n \zeta(\{ 2k \}_n).
\end{align}
In general,
\begin{align}\label{eq:zetamkn}
\sum_{\substack{p_1+\cdots+p_m=mn, \\ p_1,\ldots,p_m \geq 0}}
                          \omega^{0\cdot p_1 + 1\cdot p_2+\cdots+(m-1)\cdot p_m}
              \zeta(\{k\}_{p_1})\cdots\zeta(\{k\}_{p_m}) = (-1)^{mn-n} \zeta(\{mk\}_n).
\end{align}
\end{cor}

\section{Variation of Zagier-Broadhurst's formula}

\subsection{The shuffle automaton of $(-xy)^{\star} \sh (xy)^{\star}x^2(xy)^{\star}$}

As an application of Waldschmidt's formula \eqref{Waldschmidt}, we show 

\begin{thm}
The following relations of MZVs hold:
\begin{enumerate}[$(i)$]
\item 
$\displaystyle \hspace{8mm}
\pi^{4n+2}\sum_{p=0}^{n-1}(-1)^{p+1}\frac{4^{2n-p+1}B_{4n-4p+2}}{(4p+2)!(4n-4p+2)!} 
                               \hspace{20mm} \text{(weight$=4n+2$)}$
\begin{align}\label{thm:Sawada1}
   =-2\sum_{p=0}^{n-1}(-4)^{p}\zeta (\{ 3,1 \}_p,3,3,\{2\}_{2(n-p-1)})
                             &- 3\sum_{p=0}^{n-1}(-4)^{p}\zeta (\{3,1\}_p,4,\{2\}_{2(n-p)-1})
                                                              \quad\quad  \\
&+2(-4)^n\sum_{p=0}^{n-1}\zeta (\{3,1\}_p,5,1,\{3,1\}_{n-p-1}). \nonumber
\end{align}
\item 
$\displaystyle \hspace{8mm}
\frac{(-1)^{n}2\cdot 4^{n+1}(n+1)\pi^{4n+4}}{(4n+6)!}  
                              \hspace{20mm}   \text{(weight$=4n+4$)}$
\begin{align}\label{thm:Sawada2}
   & =  2\sum_{p=0}^{n-1}(-4)^p\zeta(\{ 3,1 \}_p,3,3,\{2\}_{2(n-p)-1})
                +3\sum_{p=0}^{n}(-4)^p\zeta (\{3,1\}_p,4,\{2\}_{2(n-p)}) \\
   &-2(-4)^{n}\sum_{p=0}^{n-1}\zeta (\{3,1\}_p,3,4,1,\{3,1\}_{n-p-1})
                  -2(-4)^{n}\sum_{p=0}^{n}\zeta (\{3,1\}_p,4,\{3,1\}_{n-p}). \nonumber
\end{align}
\end{enumerate}
where $B_k$'s are the Bernoulli numbers;
\begin{align}\label{def:Bernoulli}
                \sum_{k=0}^{\infty} B_k \frac{t^k}{k!} 
                     = \frac{te^t}{e^t-1}.
\end{align}
\end{thm}

These relations are derived from the following shuffle automata:

\begin{prop}
We have
\begin{enumerate}[$(i)$]
\item
$\displaystyle \hspace{8mm}
      (-xy)^{\star}\sh(xy)^{\star}x^2 $
\begin{align}\label{sh:sawada1}
     = (-4x^2y^2)^{\star}x^2(-xy)^{\star}-2(-4x^2y^2)^{\star}x^2yx(-xy)^{\star} 
         -2(-4x^2y^2)^{\star}x^3y(-xy)^{\star},  
\end{align}
\item
$\displaystyle \hspace{8mm}
(-xy)^{\star}\sh (xy)^{\star}x^2(xy)^{\star}$
\begin{align}\label{sh:sawada2}
&= 
(-4x^2y^2)^{\star}x^2(-4x^2y^2)^{\star} - 4(-4x^2y^2)^{\star}x^4y^2(-4x^2y^2)^{\star} 
                               \hspace{35mm}  \\
& \hspace{15mm} -2(-4x^2y^2)^{\star}x^2yx(-4x^2y^2)^{\star}
                        -2(-4x^2y^2)^{\star}x^3y(-4x^2y^2)^{\star}.
                                                            \nonumber
\end{align}
\end{enumerate}
\end{prop}
\begin{pf}
$(i)$ \quad The transition diagram of the shuffle automaton $(-xy)^{\star}\sh(xy)^{\star}x^2$ is
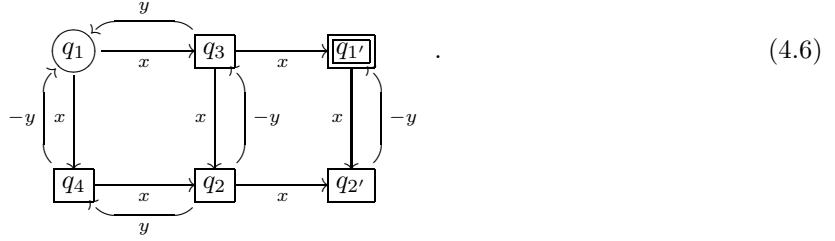
\begin{align}
\UseTips
\xymatrix @+4mm {
 *++[o][F]{q_1} \ar[r]_{x} \ar[d]_{x}
& *+[F]{q_3} 
\ar `ul^l[l]+/u4mm/`l^dl[l]_{y} [l]  
\ar[d]_{x} 
\ar[r]_{x}
& *+[F=]{q_{1'}}
\ar[d]_{x}
\\
 *+[F]{q_4} 
 \ar[r]_{x}
 \ar `ul_u[u]+/l4mm/`u_ru[u]^{-y} [u] 
& *+[F]{q_2} 
 \ar `dl_l[l]+/d4mm/`l_ul[l]^{y} [l]
 \ar `ur^u[u]+/r4mm/`u^lu[u]_{-y} [u] 
 \ar[r]_{x}
& *+[F]{q_{2'}}
\ar `ur^u[u]+/r4mm/`u^lu[u]_{-y} [u] 
}\ .
\end{align}
Let $S_{j_1,j_2}^{[i_1,\ldots,i_m]}$ be the sum of words accepted by paths from
$q_{j_1}$ to $q_{j_2}$ not passing through $q_{i_1},\ldots,q_{i_m}$. Then
$S_{1,1'}$ is the element accepted by this automaton, and we have 
\begin{align*}
       S_{1,1'} \,=\, S_{1,1}\,S_{1,2}^{[1]}\,S_{2,2'}^{[3,4,1']}\,S_{2',1'}
                      \,+\,   S_{1,1}\,S_{1,3}^{[1]}\,S_{3,1'}^{[1,2,2']}\,S_{1',1'}.
\end{align*}
By Waldschmidt's formula \eqref{Waldschmidt}, $S_{1,1} = (-4x^2y^2)^{\star}$. Other terms
are calculated as
\begin{align*}
& S_{1,2}^{[1]} = 2x^2(1+yxS_{2,2}^{[1]}-yxS_{2,2}^{[1]}) = 2x^2,   \hspace{30mm} \\
& S_{2,2'}^{[3,4,1']} = x,                                         \\
& S_{2',1'} = -y(-xy)^{\star},                                       \\
& S_{1,3}^{[1]} = x + 2x^2S_{2,3}^{[1]},                             \\
& S_{2,3}^{[1]} = -y + yxS_{2,3}^{[1]} - yxS_{2,3}^{[1]} = -y,       \\
& S_{3,1'}^{[1,2,2']} = x,                                           \\
& S_{1',1'} = (-xy)^{\star}.
\end{align*}
Hence we have \eqref{sh:sawada1}. \\
$(ii)$ \quad The transition diagram of the shuffle automaton 
$(-xy)^{\star}\sh(xy)^{\star}x^2(xy)^{\star}$ is
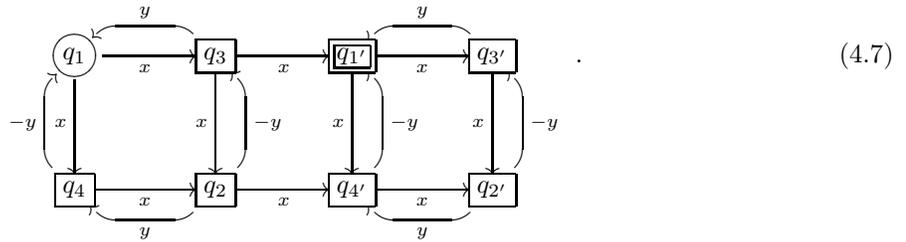
\begin{align}
\UseTips
\xymatrix @+4mm {
 *++[o][F]{q_1} \ar[r]_{x} \ar[d]_{x}
& *+[F]{q_3} 
\ar `ul^l[l]+/u4mm/`l^dl[l]_{y} [l]  
\ar[d]_{x} 
\ar[r]_{x}
& *+[F=]{q_{1'}}
\ar[r]_{x}
\ar[d]_{x}
& *+[F]{q_{3'}}
\ar[d]_{x}
\ar `ul^l[l]+/u4mm/`l^dl[l]_{y} [l]  
\\
 *+[F]{q_4} 
 \ar[r]_{x}
 \ar `ul_u[u]+/l4mm/`u_ru[u]^{-y} [u] 
& *+[F]{q_2} 
 \ar `dl_l[l]+/d4mm/`l_ul[l]^{y} [l]
 \ar `ur^u[u]+/r4mm/`u^lu[u]_{-y} [u] 
 \ar[r]_{x}
& *+[F]{q_{4'}}
\ar[r]_{x}
\ar `ur^u[u]+/r4mm/`u^lu[u]_{-y} [u] 
& *+[F]{q_{2'}} 
\ar `ur^u[u]+/r4mm/`u^lu[u]_{-y} [u] 
\ar `dl_l[l]+/d4mm/`l_ul[l]^{y} [l]
 }\ .
\end{align}
By similar consideration, we obtain \eqref{sh:sawada2}.
\hfill$\square$
\end{pf}

\newtheorem{sawadapf}{Proof of Theorem 4.1}
\renewcommand{\thesawadapf}{}

\begin{sawadapf} 
Let $(w)^+ = w(w)^{\star}$ for a word $w$.
From \eqref{sh:sawada1} and \eqref{sh:sawada2}, we have
\begin{align*}
    & (-xy)^{\star} \sh (xy)^{\star}x^2(xy)^+  \hspace{50mm} \\
    & \hspace{20mm} =  (-xy)^{\star} \sh (xy)^{\star}x^2(xy)^{\star} -  
                               (-xy)^{\star} \sh (xy)^{\star}x^2 \\
    & \hspace{20mm} = (-4x^2y^2)^{\star }x^2(-4x^2y^2)^{+ } - 
                                  4(-4x^2y^2)^{\star }x^4y^2(-4x^2y^2)^{\star }\\
    & \hspace{20mm} -2(-4x^2y^2)^{\star }x^2yx(-4x^2y^2)^{+} - 
                                             2(-4x^2y^2)^{\star }x^3y(-4x^2y^2)^{\star }\\
    & \hspace{20mm} +2(-4x^2y^2)^{\star}x^3y(-xy)^{\star} + 
                         2(-4x^2y^2)^{\star}x^2yx(-xy)^{+}-(-4x^2y^2)^{\star}x^2(-xy)^{+}.
\end{align*}
Note that, since the shuffle product preserves the length of words, 
the elements of the same length in both sides are equal.
Picking up the terms of length $4n+2$ from the RHS, we have
\begin{align*}
    & Z( \text{the sum of the terms of length} \ 4n+2 \ \text{in the RHS})  \hspace{20mm} \\
    & = -2\sum_{p=0}^{n-1}(-4)^{p}\zeta (\{ 3,1 \}_p,3,3,\{2\}_{2(n-p-1)})
                 -3\sum_{p=0}^{n-1}(-4)^{p}\zeta (\{3,1\}_p,4,\{2\}_{2(n-p)-1})\\
      & \hspace{50mm} +2(-4)^n\sum_{p=0}^{n-1}\zeta (\{3,1\}_p,5,1,\{3,1\}_{n-p-1}).
\end{align*}
On the other hand, by \eqref{eq:homsh}
	\begin{align*}
     & Z( \text{the sum of the terms of length} \ 4n+2 \ \text{in the LHS})  \hspace{20mm} \\
     & \hspace{10mm} = \sum_{p=1}^{2n}(-1)^{p-1}\left( \zeta (\{2\}_{p-1})
                         \sum_{q=1}^{2n-p+1}\zeta(\{2\}_{q-1},4,\{2\}_{2n-p+1-q}) \right).
\end{align*}

\begin{lem}
We have
\begin{align}\label{lem:sawada1}
\sum_{p=1}^{2n}(-1)^{p-1}\left( \zeta (\{2\}_{p-1})
                         \sum_{q=1}^{2n-p+1}\zeta(\{2\}_{q-1},4,\{2\}_{2n-p+1-q}) \right)
                         = (-1)^n\sum_{p=1}^{n}\zeta (\{4\}_{p-1},6,\{4\}_{n-p})
\end{align}
\end{lem}

\begin{pf}
By the definition of the harmonic product \eqref{def:harmonic}, we have
\begin{align*}
    \sum_{q=1}^{2n-p+1}\zeta(\{2\}_{q-1},4,\{2\}_{2n-p+1-q})
                         = \zeta(2)\zeta(\{2\}_{2n-p+1})-(2n-p+2)\zeta(\{2\}_{2n-p+2}).
\end{align*}
Hence
\begin{align*}
      \text{the LHS of \eqref{lem:sawada1}} = \sum_{p=1}^{2n}(-1)^{p-1}
      \left\{ \zeta(2)\zeta(\{2\}_{p-1})\zeta(\{2\}_{2n-p+1}) 
      - (2n-p+2)\zeta(\{2\}_{p-1})\zeta(\{2\}_{2n-p+2}) \right\}. 
\end{align*}
By using \eqref{eq:zeta2kn} and \eqref{def:harmonic}, it is easy to see
\begin{align*}
        & \text{the LHS of \eqref{lem:sawada1}}=  
         (-1)^n \sum_{p=1}^n \zeta(\{4\}_{p-1},6,\{4\}_{n-p})+
          (-1)^n \sum_{p=1}^{n+1}\zeta(\{4\}_{p-1},2,\{4\}_{m+1-p}) \\
        & \hspace{40mm}    + \sum_{p=0}^n(-1)^{p+1}(2n+1-2p)\zeta(\{2\}_p)\zeta(\{2\}_{2n+1-p}).
\end{align*}
Furthermore, by induction, one can prove 
\begin{align}
     (-1)^{n+1} \sum_{p=1}^{n+1} z_4^{p-1} z_2 z_4^{n+1-p} =
                \sum_{p=0}^n (-1)^{p+1} (2n+1-2p) z_2^p * z_2^{2n+1-p}.
\end{align}
Thus the proof is completed.
\hfill$\square$
\end{pf}

By calculating harmonic product, we have
\begin{align*}
& (-1)^n\sum_{p=1}^{n}\zeta (\{4\}_{p-1},6,\{4\}_{n-p})
=(-1)^n\left(\zeta(6)\zeta (\{4\}_{n-1}) - 
                         \sum_{p=1}^{n-1}\zeta (\{4\}_{p-1},10,\{4\}_{n-p-1})\right)\\
& =(-1)^n\zeta(6)\zeta (\{4\}_{n-1}) + 
        (-1)^{n+1}\zeta(10)\zeta (\{4\}_{n-2})+
                          (-1)^{n+2}\sum_{p=1}^{n-2}\zeta (\{4\}_{p-1},10,\{4\}_{n-p-2})\\
& \vdots \\
& =\sum_{p=0}^{n-1}(-1)^{p+1}\zeta (\{4\}_{p})\zeta (4n-4p+2).
\end{align*}
Here substituting the formulas 
\begin{align}\label{eq:zeta2k}
\zeta(2k) = \frac{(-1)^{k-1}(2\pi)^{2k}B_{2k}}{2(2k)!},
\end{align}
\begin{align}\label{eq:zeta4kn} 
\zeta (\{4\}_{n})=\frac{2 \cdot 4^n{\pi}^{4n}}{(4n+2)!},  
\end{align}
we obtain the LHS of \eqref{thm:Sawada1}. 

We can prove\ \eqref{thm:Sawada2}
in a similar way. This completes the proof.
\hfill$\square$
\end{sawadapf}

\subsection{The shuffle automaton of $(x^2y)^{\star} \sh (-x^2y)^{\star}$}

Through consideration on the shuffle automaton of
$(x^2y)^{\star} \sh (-x^2y)^{\star}$, we show

\begin{thm}\label{thm:x2y}
We have
\begin{align}\label{thm:kitani1}
\sum_{\varepsilon_i, \, \varepsilon'_j=0,1}
   \frac{12^n}{2^{\varepsilon_1+\cdots+\varepsilon_n+\varepsilon'_1+\cdots+\varepsilon'_{n-1}}}
   \zeta(\{ 5-\varepsilon'_{i-1}-\varepsilon_i ,1+\varepsilon_i+\varepsilon'_i \}_{i=1}^n ) 
  = \frac{6(2\pi)^{6n}}{(6n+3)!}
\end{align}
where $\varepsilon'_0=\varepsilon'_n=0$.
\end{thm}

First we show the following lemma:

\begin{lem}\label{lem:wsh}
Assume that $w \in \mathfrak{H}$ be a word of length $l$. Let 
$A=(a_{ij})$  be the adjacency matrix of the finite automaton of 
$w^{\star} \sh (-w)^{\star}$ and $A^n=(a_{ij}^{(n)})$. Then
\begin{align*}
      n \notin 2l \mathbb{Z}_{\ge 0} \  \Longrightarrow  \ a_{11}^{(n)}=0.
\end{align*}
\end{lem}
\begin{pf}
It is obvious that $n \notin l \mathbb{Z}_{\ge 0} \ \Longrightarrow  \ a_{11}^{(n)}=0$.
Suppose $n=(2k+1)l \ (k \in \mathbb{Z}_{\ge 0})$. Then
\begin{align*}
\sum_{i=0}^{2k+1} w^i \sh (-w)^{2k+1-i} 
&= \sum_{i=0}^{k} \{ w^i \sh (-w)^{2k+1-i} + w^{2k+1-i} \sh (-w)^i \}\\
&= \sum_{i=0}^{k} \{ (-1)^{2k+1-i}+(-1)^i\} \; w^i \sh w^{2k+1-i} \\
&= 0.
\end{align*}
\hfill$\square$
\end{pf}

The shuffle automaton of $(x^2y)^{\star} \sh (-x^2y)^{\star}$ is represented as

\begin{align}\label{fig:w33}
\UseTips
\xymatrix @+6mm {
 *++[o][F=]{q_1}
 \ar `ur_r[r]+/u4mm/`r_dr[r]^{x} [r]
 \ar `dl^d[d]+/l4mm/`d^dr[d]_{x} [d] 
& *+[F]{q_4}
 \ar `ur_r[r]+/u4mm/`r_dr[r]^{x} [r]
\ar `dr_d[d]+/r4mm/`d_dl[d]^{x} [d] 
& *+[F]{q_7} 
\ar `u^l[ll]+/u1cm/`l^d[ll]_{y} [ll]  
\ar `dr_d[d]+/r4mm/`d_dl[d]^{x} [d] 
\\
 *+[F]{q_5}
 \ar `dr^r[r]+/u4mm/`r^ur[r]_{x} [r]
 \ar `dl^d[d]+/l4mm/`d^dr[d]_{x} [d] 
& *+[F]{q_8}
\ar `dr^r[r]+/u4mm/`r^ur[r]_{x} [r]
  \ar `dr_d[d]+/r4mm/`d_dl[d]^{x} [d] 
& *+[F]{q_2} 
\ar `ul^l[ll]+/u6mm/`l^dl[ll]_{y} [ll]  
\ar `dr_d[d]+/r4mm/`d_dl[d]^{x} [d] 
\\
 *+[F]{q_9} 
  \ar `dr^r[r]+/u4mm/`r^ur[r]_{x} [r]
 \ar `l_u[uu]+/l1cm/`u_r[uu]^{-y} [uu] 
& *+[F]{q_3}
  \ar `dr^r[r]+/u4mm/`r^ur[r]_{x} [r]
 \ar `lu_u[uu]+/l6mm/`u_ru[uu]^{-y} [uu]
& *+[F]{q_6} 
 \ar `d_l[ll]+/d1cm/`l_u[ll]^{y} [ll]
 \ar `r^u[uu]+/r1cm/`u^l[uu]_{-y} [uu] 
}
\end{align}

The adjacency matrix is
\begin{align}\label{eq:Ashape1}
A= 
\begin{bmatrix}
0 & P_1 & 0 \\
0 & 0 & P_2 \\
P_3 & 0 & 0 
\end{bmatrix}
\end{align}
where
\begin{align*}
P_1=
\begin{bmatrix}
x & x & 0 \\
0 & y & x \\
-y & 0 & x 
\end{bmatrix}
,P_2=
\begin{bmatrix}
x & x & 0 \\
0 & x & x \\
-y & 0 & y 
\end{bmatrix}
,P_3=
\begin{bmatrix}
y & x & 0 \\
0 & x & x \\
-y & 0 & x 
\end{bmatrix}.
\end{align*}
From Lemma \ref{lem:wsh} and Propostion \ref{prop:am}, the words accepted by this automaton are
$a_{11}^{(6n)}\; (n \in \mathbb{Z}_{\ge 0})$. Because of the form of $A$ \eqref{eq:Ashape1},
it is the $(1,1)$-entry of $\{P_1P_2P_3\}^{2n}$.
Compute $P:=\{P_1P_2P_3\}^2$:

\begin{align}\label{eq:P}
P
&=
\left\{
\begin{bmatrix}
x & x & 0 \\
0 & y & x \\
-y & 0 & x 
\end{bmatrix}
\begin{bmatrix}
x & x & 0 \\
0 & x & x \\
-y & 0 & y 
\end{bmatrix}
\nonumber
\begin{bmatrix}
y & x & 0 \\
0 & x & x \\
-y & 0 & x 
\end{bmatrix} 
\right\}^2 \nonumber \\
%&=\begin{bmatrix}
%0 & 3x^3 & 3x^3 \\
%-(2xy^2+yxy) & -(xyx-yx^2) & xyx+2yx^2 \\
% -(2xy^2+yxy) & -(xyx+2yx^2) & xyx-yx^2
%\end{bmatrix}^2  \nonumber \\
&= \begin{bmatrix}
p & q & -q \\
r & s & t \\
-r & t & s
\end{bmatrix}
\end{align}
where
\begin{align}
\left\{
\begin{array}{rl}
p&=-3(4x^4y^2+2x^3yxy), \\
q&=-3(2x^4y+x^3yx^2), \\
r&=-3(2yx^3y^3+yx^2yxy), \\
s&=-3(2xy^2x^3+yxyx^3+yx^3yx+xyxyx^2+yx^2yx^2), \\
t&=-3(2xy^2x^3+yxyx^3-yx^3yx+xyxyx^2). 
\end{array}\right.
\end{align}

\begin{lem}\label{lem:w3}
There exist elements $p_n, q_n, r_n, s_n, t_n \in \mathfrak{H}$
such that 
\begin{align*}
   P^n = 
   \begin{bmatrix}
   p_n  & q_n & -q_n \\
   r_n  & s_n & t_n  \\
   -r_n & t_n & s_n
\end{bmatrix}.
\end{align*}
\end{lem}
\begin{pf}
From \eqref{eq:P}, this statement is true for $n=1$. Now compute $P^2$: 
\begin{align*}
P^2=\begin{bmatrix}
p^2+2qr & pq+q(s-t) & -\{pq+q(s-t) \} \\
rp+(s-t)r & rq+s^2+t^2 & -rq+st+ts \\
-\{rp+(s-t)r\} & -rq+st+ts & rq+s^2+t^2 
\end{bmatrix}.
\end{align*}
Hence it is represented as
\begin{align*}
\begin{bmatrix}
p_2 & q_2 & -q_2 \\
r_2 & s_2 & t_2 \\
-r_2 & t_2 & s_2
\end{bmatrix}.
\end{align*}
For $n \ge 3$, it is proved by induction.
\hfill$\square$
\end{pf}

From this lemma and $P^n=P^{n-1}P$, we have the following recursive equations
for $p_n, \ q_n$:
\begin{align}\label{eq:zenkashikiPn}
\left\{
\begin{array}{rl}
p_n &= -6p_{n-1}(2x^4y^2+x^3yxy) - 6q_{n-1}y(2x^3y^2+x^2yxy), \\
   \\
q_n &= -3p_{n-1}(2x^4yx+x^3yx^2)-3q_{n-1}y(2x^3yx+x^2yx^2).
\end{array}
\right. 
\end{align}
Set $w_n=(x^4y^2)^n$ and $w_n'=(x^4y^2)^{n-1}x^4yx \ (w_0=w'_0=1)$.
We define the operations $\sigma_i,\tau_j$ by
\begin{align}
\left\{\begin{array}{rl}
\sigma_i(w_n^{(')}) 
&=w_{i-1} \overbrace{x^3\underline{yx}y}^{i} w_{n-i}^{(')} \quad(1 \leq i \leq n-1), \quad
      \sigma_n(w_n')=w_{n-1}x^3\underline{yx}x \\
\\
\tau_j(w_n^{(')})
&=w_{j-1} \overbrace{x^4 y\underline{xy}x^3 y^2}^{j,j+1} w_{n-1-j}^{(')}
\quad (1 \leq j \leq n-2) \quad 
\tau_{n-1}(w_n')=w_{n-2}x^4y\underline{xy}yx.
\end{array}
\right.
\end{align}
where the underlines designate the position of exchanging the order of $x$ and $y$.
Solving the recursive equatoins \eqref{eq:zenkashikiPn} in terms of these operations, we have the 
following proposition:

\begin{prop}
\begin{align}\label{eq:pnqn}
p_n&= \sum_{\substack{\varepsilon_1,\ldots,\varepsilon_n= 0,1\\
\varepsilon_1^{\prime},\ldots,\varepsilon_{n-1}^{\prime}=0,1}}
\frac{(-12)^n}{2^{\varepsilon_1+\cdots+\varepsilon_n+\varepsilon'_1
+\cdots+\varepsilon'_{n-1}}}
\sigma_1^{\varepsilon_1}\cdots\sigma_n^{\varepsilon_n}
\tau_1^{\varepsilon'_1}\cdots\tau_{n-1}^{\varepsilon'_{n-1}}
(w_n),\\
q_n&= \sum_{\substack{\varepsilon_1,\ldots,\varepsilon_n= 0,1\\
\varepsilon_1^{\prime},\ldots,\varepsilon_{n-1}^{\prime}=0,1}}
\frac{(-12)^n}{2^{1+\varepsilon_1+\cdots+\varepsilon_n
+\varepsilon'_1+\cdots+\varepsilon'_{n-1}}}
\sigma_1^{\varepsilon_1}\cdots\sigma_n^{\varepsilon_n}
\tau_1^{\varepsilon'_1}\cdots\tau_{n-1}^{\varepsilon'_{n-1}}
(w_n').
\end{align}
\end{prop}

\begin{pf}
We prove by induction. For $n=1$, we have
\begin{align*}
p_1 &= -12x^4y^2-6x^3yxy=-12\left( x^4y^2+\frac{\sigma_1}{2}x^4y^2 \right), \\
q_1 &= -6x^4yx-3x^3yx^2=-12\left( \frac{1}{2}x^4yx+\frac{\sigma_1}{4}x^4yx \right).
\end{align*}
So the statement is true. Next assume \eqref{eq:pnqn} to be true for $p_n, \ q_n$. 
As for $p_{n+1}$, from \eqref{eq:pnqn}, it follows that
\begin{align*}
p_{n+1}&= -12p_nx^4y^2-6p_nx^3yxy - 12q_nyx^3y^2 -6q_nyx^2yxy \\
&= \sum_{\substack{\varepsilon_1,\ldots,\varepsilon_n= 0,1\\
\varepsilon_1^{\prime},\ldots,\varepsilon_{n-1}^{\prime}=0,1}}
\frac{(-12)^{n+1}}{2^{\varepsilon_1+\cdots+\varepsilon_n
+\varepsilon'_1+\cdots+\varepsilon'_{n-1}}}
\sigma_1^{\varepsilon_1}\cdots\sigma_n^{\varepsilon_n}
\tau_1^{\varepsilon'_1}\cdots\tau_{n-1}^{\varepsilon'_{n-1}} (x^4y^2)^{n+1} \\
&+ \sum_{\substack{\varepsilon_1,\ldots,\varepsilon_n= 0,1\\
\varepsilon_1^{\prime},\ldots,\varepsilon_{n-1}^{\prime}=0,1}} 
\frac{(-12)^{n+1}}{2^{1+\varepsilon_1+\cdots+\varepsilon_n
+\varepsilon'_1+\cdots+\varepsilon'_{n-1}}}
\sigma_1^{\varepsilon_1}\cdots\sigma_n^{\varepsilon_n}
\tau_1^{\varepsilon'_1}\cdots\tau_{n-1}^{\varepsilon'_{n-1}}
\sigma_{n+1} (x^4y^2)^{n+1} \\
&+ \sum_{\substack{\varepsilon_1,\ldots,\varepsilon_n= 0,1 \\
\varepsilon_1^{\prime},\ldots,\varepsilon_{n-1}^{\prime}=0,1}}
\frac{(-12)^{n+1}}{2^{1+\varepsilon_1+\cdots+\varepsilon_n
+\varepsilon'_1+\cdots+\varepsilon'_{n-1}}}
\sigma_1^{\varepsilon_1}\cdots\sigma_n^{\varepsilon_n}
\tau_1^{\varepsilon'_1}\cdots\tau_{n-1}^{\varepsilon'_{n-1}}
\tau_n (x^4y^2)^{n+1} \\
&+ \sum_{\substack{\varepsilon_1,\ldots,\varepsilon_n= 0,1 \\
\varepsilon_1^{\prime},\ldots,\varepsilon_{n-1}^{\prime}=0,1}}
\frac{(-12)^{n+1}}{2^{2+\varepsilon_1+\cdots+\varepsilon_n
+\varepsilon'_1+\cdots+\varepsilon'_{n-1}}}
\sigma_1^{\varepsilon_1}\cdots\sigma_n^{\varepsilon_n}
\tau_1^{\varepsilon'_1}\cdots\tau_{n-1}^{\varepsilon'_{n-1}}
\sigma_{n+1}\tau_n (x^4y^2)^{n+1} \\
&=\sum_{\substack{\varepsilon_1,\ldots,\varepsilon_{n+1}= 0,1 \\
\varepsilon_1^{\prime},\ldots,\varepsilon_n^{\prime}=0,1}}
\frac{(-12)^{n+1}}{2^{\varepsilon_1+\cdots+\varepsilon_{n+1}
+\varepsilon'_1+\cdots+\varepsilon'_n}}
\sigma_1^{\varepsilon_1}\cdots\sigma_{n+1}^{\varepsilon_{n+1}}
\tau_1^{\varepsilon'_1}\cdots\tau_n^{\varepsilon'_n}
(w_{n+1}).
\end{align*}
Hence \eqref{eq:pnqn} holds for $p_{n+1}$.
The case for $q_{n+1}$ is proved in a similar way.
\hfill $\square$
\end{pf}

\newtheorem{kitanipf}{Proof of Theorem \ref{thm:x2y}}
\renewcommand{\thekitanipf}{}

\begin{kitanipf}
Transforming the RHS of \eqref{eq:pnqn} via $Z$, we obtain the LHS \eqref{thm:kitani1}.
By the finite double shuffle relation Proposition \ref{prop:fdsr} and 
the formula \eqref{eq:zeta2kn} (Corollary \ref{cor:zhkk}), we have
\begin{align*}
    Z(\text{the sum of the terms of} \ (x^2y)^{\star} \sh (-x^2y)^{\star} \text{of length 6n}\, )
    = \zeta(\{6\}_n).
\end{align*}
Furthermore, we have the formula (see Appendix)
\begin{align}\label{eq:zeta6n}
         \zeta(\{6\}_n)= \frac{6(2\pi)^{6n}}{(6n+3)!} \ .
\end{align}
Thus the proof is completed.
\hfill$\square$
\end{kitanipf}

\subsection{The shuffle automaton of $(xy)^{\star} \sh (\omega xy)^{\star} \sh (\omega^2xy)^{\star}$}

Let $\omega$ be a primitive cubic root of unity. Through considration on the shuffle automaton of
$(xy)^{\star} \sh (\omega xy)^{\star} \sh (\omega^2xy)^{\star}$,

\begin{thm}\label{thm:dim3}
We have
\begin{align}\label{eq:dim3}
\sum_{\substack{\varepsilon_1,\ldots,\varepsilon_n= 0,1\\
\varepsilon_1^{\prime},\ldots,\varepsilon_{n-1}^{\prime}=0,1}}
&\frac{36^n}{3^{\varepsilon_1+\cdots+\varepsilon_n
+\varepsilon'_1+\cdots+\varepsilon'_{n-1}}}
\zeta( \{4-\varepsilon_i-\varepsilon_{i-1}',1+\varepsilon_i,1 +\varepsilon_i' \}_{i=1}^n)
=\frac{6(2\pi)^{6n}}{(6n+3)!} \ .
\end{align}
where $\varepsilon_0'=\varepsilon_n'=0$.
\end{thm}

The transition diagram and the adjacency matrix of this shuffle automaton are as follows:

\begin{align}\label{fig:dim3}
\UseTips
\xymatrix @+8mm {
 *++[o][F=]{q_1} \ar[rrr]_{x} \ar[dr]^{x}
\ar `l^d[ddd]+/l18mm/`d^r[ddd]_{x} [ddd]  
& & 
& *+[F]{q_6} 
\ar `ul^l[lll]+/u6mm/`l^dl[lll]_{\omega y} [lll]  
\ar `r_d[ddd]+/r18mm/`d_l[ddd]^{x} [ddd]  
\ar[dl]_{x} 
\\
& *+[F]{q_7} \ar[r]_{x} \ar[d]^{x}
 \ar `l_ul[ul]+/dl5mm/`ul_u[ul]^{\omega^2 y}  [ul]  
& *+[F]{q_4} \ar `ul^l[l]+/u4mm/`l^dl[l]_{\omega y} [l]  
 \ar `r^ur[ru]+/dr5mm/`ru^u[ru]_{\omega^2 y}  [ru]  
\ar[d]_{x} 
& \\
& *+[F]{q_3} \ar[r]^{x}
 \ar `ul_u[u]+/l4mm/`u_ru[u]^{y} [u] 
 \ar `l^dl[dl]+/lu5mm/`dl^d[dl]_{\omega^2 y}  [dl]  
& *+[F]{q_8}  
 \ar `dl_l[l]+/d4mm/`l_ul[l]^{\omega y} [l]
 \ar `ur^u[u]+/r4mm/`u^lu[u]_{y} [u] 
\ar `r_dr[dr]+/ru5mm/ `dr_d [dr]^{\omega^2 y}  [dr]  
& \\
 *+[F]{q_5}
 \ar[rrr]^{x} \ar[ur]_{x} 
 \ar `ul_u[uuu]+/l6mm/`u_ru[uuu]^{y} [uuu] 
& & 
& *+[F]{q_2} 
 \ar `dl_l[lll]+/d6mm/`l_ul[lll]^{\omega y} [lll]
 \ar `ur^u[uuu]+/r6mm/`u^lu[uuu]_{y} [uuu] 
\ar[ul]^{x} 
}
\end{align}

\begin{align}\label{matrix:dim3}
A=\begin{bmatrix}
0&0&0&0&x&x&x&0 \\
0&0&0&0&\omega y&y&0&x \\
0&0&0&0&\omega^2 y&0&y&x \\
0&0&0&0&0&\omega^2 y&\omega y&x \\
y&x&x&0&0&0&0&0 \\
\omega y&x&0&x&0&0&0&0 \\
\omega^2 y&0&x&x&0&0&0&0 \\
0&\omega^2 y&\omega y&y&0&0&0&0
\end{bmatrix}.
\end{align}
Let $A^n=(a_{ij}^{(n)})$. What we want to know is $a_{11}^{(n)}$. The following lemma is
an analogy of Lemma \ref{lem:wsh}.

\begin{lem}
Let $w \in \mathfrak{H}$ be a word of length $l$. Let $A=(a_{ij})$ be the adjacency matrix
of the shuffle automaton of $w^{\star} \sh (\omega w)^{\star} \sh (\omega^2 w)^{\star}$ and
$A^n =(a_{ij}^{(n)})$. Then
\begin{align*}
   n \notin 3l\mathbb{Z} \Longrightarrow a_{11}^{(n)}=0.
\end{align*}
\end{lem}
%\begin{pf}
%\begin{align}\label{lem:dim3}
%\sum_{i+j+k=m, \, i,j,k \geq 0} \omega^{0i+1j+2k} =
%\begin{cases}
%1 & m \in 3\mathbb{Z}_{\ge 0}  \\ 
%0 & m \notin 3\mathbb{Z}_{\ge 0} 
%\end{cases}.
%\end{align}
%\begin{align*}
%\sum_{i=0}^{\infty} t^{il} \sum_{j=0}^{\infty} \omega^j t^{jl}
%\sum_{k=0}^{\infty} \omega^{2k} t^{kl}
%&=\frac{1}{1-t^l}\cdot \frac{1}{1-\omega t^l}\cdot \frac{1}{1-\omega^2 t^l} \\
%&=\frac{1}{1-t^{3l}}   \\
%&=\sum_{m=0}^{\infty} t^{3lm}.
%\end{align*}
%On the other hand, 
%\begin{align*}
% \sum_{i=0}^{\infty} t^{il} \sum_{j=0}^{\infty} \omega^j t^{jl}
%\sum_{k=0}^{\infty} \omega^{2k} t^{kl}
%=\sum_{m=0}^{\infty} \left( \sum_{i+j+k=m} \omega^{0i+1j+2k} \right) t^{lm}. 
%\end{align*}
%Compare these, we obtain the \eqref{lem:dim3}.
%\hfill$\square$
%\end{pf}

\begin{pf}
Let us show $a_{11}^{(l(3m+1))}=0$:
\begin{align*}
a_{11}^{(l(3m+1))}
& =\sum_{i_1+i_2+i_3=3m+1}
w^{i_1}\sh(\omega w)^{i_2} \sh(\omega^2 w)^{i_3} \\
&= \sum_{i_1+i_2+i_3=3m+1}\omega^{i_2+2i_3} 
w^{i_1}\sh w^{i_2} \sh w^{i_3} \\
&=\sum_{\substack{j_1+j_2+j_3=3m+1 \\ j_1>j_2>j_3}}
\alpha(\omega) \; w^{j_1}\sh w^{j_2} \sh w^{j_3} 
+\sum_{\substack{2j_1+j_2=3m+1 \\ j_1 \neq j_2}}
\beta(\omega) \; w^{j_1}\sh w^{j_1} \sh w^{j_2}, 
\end{align*}
where
\begin{align*}
\left\{ \begin{array}{rl}
\alpha(\omega)&=\omega^{2j_1+j_2} +\omega^{2j_1+j_3}+ \omega^{2j_2+j_1} +
\omega^{2j_2+j_3} +\omega^{2j_3+j_1} +\omega^{2j_3+j_2}, \\
                                                         \\
\beta(\omega)&=\omega^{j_1+2j_2} +\omega^{j_1+2j_1}+ \omega^{j_2+2j_1}
\end{array} \right.
\end{align*}
%Using $\omega^3=1$ and $\omega^2 + \omega + 1=0$, we have
%\begin{align*}
%\alpha(\omega)&= \omega^{3m} \omega^{2j_1+j_2} +\omega^{3m} \omega^{2j_1+j_3}
%+ \omega^{6m} \omega^{2j_2+j_1} + \omega^{2j_2+j_3} 
%+ \omega^{6m}\omega^{2j_3+j_1} +\omega^{2j_3+j_2}  \\
%&= \omega^{ (j_1+j_2+j_3-1) +2j_2+j_1} + \omega^{(j_1+j_2+j_3-1) +2j_1+j_3}
%+ \omega^{2(j_1+j_2+j_3-1) +2j_2+j_1} \\
%&\hspace{4.9cm} + \omega^{2j_2+j_3} 
%+\omega^{2(j_1+j_2+j_3-1) +2j_3+j_1} +\omega^{2j_3+j_2}  \\
%&= \omega^{ 2j_2+j_3-1} + \omega^{j_2+2j_3-1}
%+ \omega^{j_2+2j_3-2} + \omega^{2j_2+j_3} 
%+\omega^{2j_2+j_3-2} +\omega^{2j_3+j_2}  \\
%&= \omega^{2j_2+j_3-3} (1+\omega +\omega^2 )
%+\omega^{j_2+2j_3-3} (1+\omega +\omega^2 )\\
%&= 0 \\
%\beta(\omega)&= \omega^{3m}\omega^{j_1+2j_2} +1+ \omega^{6m}\omega^{j_2+2j_1} \\
%&=1+\omega^{ (2j_1+j_2-1) +j_1+2j_2}+ \omega^{2(2j_1+j_2-1) +j_2+2j_1} \\
%&=1+\omega^{-1}+ \omega^{-2} \\
%&=0
%\end{align*}
One can easily show $\alpha(\omega)=\beta(\omega)=0$ by noting that 
$\omega^{j_1+j_2+j_3-1}=1$ for $j_1+j_2+j_3=3m+1$ and 
$\omega^{2j_1+j_2-1}=1$ for $2j_1+j_2=3m+1$.
Thus $a_{11}^{(l(3m+1))}=0$. 

One can prove $a_{11}^{(l(3m+2))}=0$ in a similar way.
\hfill$\square$
\end{pf}

From this lemma and the form of $A$ \eqref{matrix:dim3}
we see that $a_{11}^{(n)}=0 \ (n \notin 6\mathbb{Z})$ and
$a_{11}^{(6n)}$ is the $(1,1)$-entry of $\widehat{P}^n$,
\begin{align}
\widehat{P} & = 
  \left\{ 
  \begin{bmatrix}
   x&x&x&0 \\
   \omega y&y&0&x \\
   \omega^2 y&0&y&x \\
   0&\omega^2 y&\omega y&x 
   \end{bmatrix}
   \begin{bmatrix}
   y&x&x&0\\
   \omega y&x&0&x\\
   \omega^2 y&0&x&x\\
   0&\omega^2 y&\omega y&y\\
   \end{bmatrix} 
   \right\}^3  \nonumber \\
 & =  
\begin{bmatrix}
p&-\omega q&-\omega^2 q&-q \\
-\omega^2 r&s+t&s+t''&s+t' \\
-\omega r&s+t'&s+t&s+t''\\
-r&s+t''&s+t'&s+t 
\end{bmatrix}.  \label{eq:widehatP}
\end{align}
where
\begin{align}\label{eq:pqrst}
\left\{
\begin{tabular}{l}
$p=  12(3x^3y^3+x^2yxy^2) $
\hspace{1.3cm}
$ t \hspace{1.6mm} = 4yxyx^2y + 4yx^2y^2x+4yxyxyx $\\
$q=  12x^3y^2x + 4x^2yxyx$
\hspace{15.3mm}
$ t'= 4\omega yxyx^2y + 4\omega^2 yx^2y^2x $\\
$r= 12yx^2y^3 + 4yxyxy^2 $
\hspace{15.7mm}
$ t''= 4\omega^2 yxyx^2y + 4\omega yx^2y^2x  $\\
$s = 12( y^2x^2yx+xy^3yx ) + 4(y^2x^2yx+xy^2xyx+yxy^2x^2)$
\end{tabular} \right.
\end{align}
By induction, one can show 

\begin{lem}
There exist $p_n, q_n, r_n, s_n, t_n, t'_n, t''_n \in \mathfrak{H}$ such that
\begin{align*}
\widehat{P}^n =
\begin{bmatrix}
p_n&-\omega q_n&-\omega^2 q_n&-q_n \\
-\omega^2 r_n&s_n+t_n&s_n+t''_n&s_n+t'_n \\
-\omega r_n&s_n+t'_n&s_n+t_n&s_n+t''_n\\
-r_n&s_n+t''_n&s_n+t'_n&s_n+t_n 
\end{bmatrix}.
\end{align*}
The elements $p_n, q_n$ satisfy the recursive relations
\begin{align}
\left\{
\begin{array}{rl}
p_n &= 36\left\{ p_{n-1}\left( x^3y^3+\frac{1}{3}x^2yxy^2 \right) 
                         + q_{n-1} \left( yx^2y^3+\frac{1}{3}yxyxy^2 \right) \right\}, \\
                                                               \\
q_n &=  36\left\{ p_{n-1}\left( \frac{1}{3}x^3y^2x+\frac{1}{3^2}x^2yxyx \right) 
                 +  q_{n-1}\left( \frac{1}{3}yx^2y^2x+\frac{1}{3^2}yxyxyx \right) \right\}.
\end{array}
\right. 
\end{align}
\end{lem}

Set $w_n=(x^3y^3)^n$ and $w_n'=(x^3y^3)^{n-1}x^3y^2x \ (w_0=w'_0=1)$.
Let us define the operations $\sigma_i,\tau_j$ by
\begin{align}
\left\{\begin{array}{rl}
\sigma_i(w_n^{(')}) 
&=w_{i-1} \overbrace{x^2\underline{yx}y^2}^{i} w_{n-i}^{(')} \quad (1 \leq i \leq n-1), 
                     \quad \sigma_n(w_n')=w_{n-1}x^2\underline{yx}yx \\
\\
\tau_j(w_n^{(')})
&=w_{j-1} \overbrace{x^3 y^2\underline{xy}x^2 y^3}^{j,j+1} w_{n-1-j}^{(')}
\quad (1 \le j \le n-2), \quad \tau_{n-1}(w_n')=w_{n-2}x^3y^2\underline{xy}x^2y^2x .
\end{array}
\right.
\end{align}
Similarly as in the previous subsection, one can show

\begin{prop}
We have
\begin{align}
p_n&= \sum_{\substack{\varepsilon_1,\ldots,\varepsilon_n= 0,1\\
\varepsilon_1^{\prime},\ldots,\varepsilon_{n-1}^{\prime}=0,1}}
\frac{(36)^n}{3^{\varepsilon_1+\cdots+\varepsilon_n+\varepsilon'_1
+\cdots+\varepsilon'_{n-1}}}
\sigma_1^{\varepsilon_1}\cdots\sigma_n^{\varepsilon_n}
\tau_1^{\varepsilon'_1}\cdots\tau_{n-1}^{\varepsilon'_{n-1}}
(w_n), \label{eq:dim3pn} \\
q_n&= \sum_{\substack{\varepsilon_1,\ldots,\varepsilon_n= 0,1\\
\varepsilon_1^{\prime},\ldots,\varepsilon_{n-1}^{\prime}=0,1}}
\frac{(36)^n}{3^{1+\varepsilon_1+\cdots+\varepsilon_n
+\varepsilon'_1+\cdots+\varepsilon'_{n-1}}}
\sigma_1^{\varepsilon_1}\cdots\sigma_n^{\varepsilon_n}
\tau_1^{\varepsilon'_1}\cdots\tau_{n-1}^{\varepsilon'_{n-1}}
(w_n').
\end{align}
\end{prop}

\newtheorem{dimpf}{Proof of Theorem \ref{thm:dim3}}
\renewcommand{\thedimpf}{}

\begin{dimpf}
Transforming \eqref{eq:dim3pn} via $Z$ to MZVs yields the LHS of \eqref{eq:dim3}.
From Proposition \ref{prop:fdsr} and Theorem \ref{thm:hkk}, it follows that
\begin{align*}
       & Z\left( (xy)^{\star} \sh (\omega xy)^{\star} \sh (\omega^2 xy)^{\star} \right) \\
     = & Z\left( (xy)^{\star} * (\omega xy)^{\star} * (\omega^2 xy)^{\star} \right)  \\
     = & Z\left( (x^5y)^{\star} \right) \\
     = & \text{the RHS of \eqref{eq:dim3}}.
\end{align*}
\hfill$\square$
\end{dimpf}

\section{Harmonic automata}

\begin{thm}\label{thm:hauto}
The automaton accepting the harmonic product $w_1*w_2$ where
$w_1=z_{p_1} z_{p_2} \cdots z_{p_m}$ and $w_2=z_{q_1} z_{q_2} \cdots z_{q_n}$
is rpresented by the following transition diagram \eqref{fig:h}. We call this
the harmonic automaton of $w_1*w_2$:

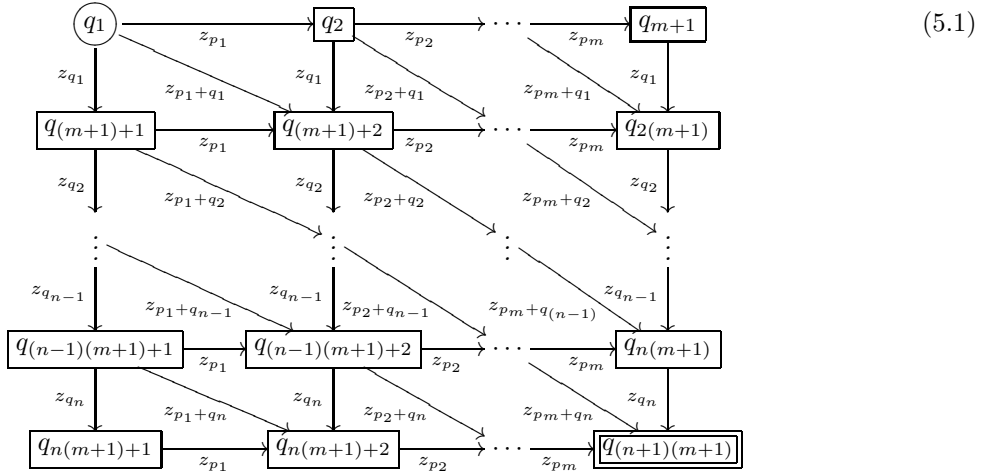
\begin{align}\label{fig:h}
\UseTips
\xymatrix @+0mm {
 *++[o][F]{q_1} \ar[r]_{z_{p_1}} \ar[d]_{z_{q_1}}@\ar[rd]_{z_{p_1+q_1}}
& *+[F]{q_2} \ar[r]_{z_{p_2}} \ar[d]_{z_{q_1}} \ar[rd]_{z_{p_2+q_1}}
& \cdots \ar[r]_{z_{p_m}} \ar[rd]_{z_{p_m+q_1}}
& *+[F]{q_{m+1}} \ar[d]_{z_{q_1}} 
\\
 *+[F]{q_{(m+1)+1}} \ar[r]_{z_{p_1}} \ar[d]_{z_{q_2}} \ar[rd]_{z_{p_1+q_2}}
& *+[F]{q_{(m+1)+2}} \ar[r]_{z_{p_2}} \ar[d]_{z_{q_2}} \ar[rd]_{z_{p_2+q_2}}
& \cdots \ar[r]_{z_{p_m}} \ar[rd]_{z_{p_m+q_2}}
& *+[F]{q_{2(m+1)}} \ar[d]_{z_{q_2}}
\\
  \vdots  \ar[d]_{z_{q_{n-1}}} \ar[rd]_{z_{p_1+q_{n-1}}}
& \vdots  \ar[d]_{z_{q_{n-1}}} \ar[rd]_{z_{p_2+q_{n-1}}}
& \vdots  \ar[rd]_{z_{p_m+q_{(n-1)}}}
& \vdots  \ar[d]_{z_{q_{n-1}}} 
\\
 *+[F]{q_{(n-1)(m+1)+1}} \ar[r]_{z_{p_1}} \ar[d]_{z_{q_n}} \ar[rd]_{z_{p_1+q_n}}
& *+[F]{q_{(n-1)(m+1)+2}} \ar[r]_{\qquad z_{p_2}} \ar[d]_{z_{q_n}} \ar[rd]_{z_{p_2+q_n}}
& \cdots \ar[r]_{z_{p_m}} \ar[rd]_{z_{p_m+q_n}}
& *+[F]{q_{n(m+1)}} \ar[d]_{z_{q_n}} 
\\
 *+[F]{q_{n(m+1)+1}} \ar[r]_{z_{p_1}} 
& *+[F]{q_{n(m+1)+2}} \ar[r]_{\quad z_{p_2}} 
& \cdots \ar[r]_{z_{p_m}\qquad} 
& *+[F=]{q_{(n+1)(m+1)}} 
}
\end{align}
\end{thm}
\begin{pf}
In the definition of harmonic product
\begin{align*}
z_{p_1} z_{p_2} \cdots z_{p_m}*z_{q_1} z_{q_2} \cdots z_{q_n}
=& z_{p_1} (z_{p_2} \cdots z_{p_m} * z_{q_1} z_{q_2} \cdots z_{q_n}) 
+ z_{q_1} (z_{p_1} z_{p_2} \cdots z_{p_m} * z_{q_2} \cdots z_{q_n}) \\
& \hspace{25mm} + z_{p_1+q_1} (z_{p_2} \cdots z_{p_m} * z_{q_2} \cdots z_{q_n}),     
\end{align*}
each term designates
\begin{align*}
&\text{(the first term)} : \text{the transition from the state} \ q_1 \ \text{to the state} \ q_2 \
                           \text{inputing} \  z_{p_1}, \\               
&\text{(the second term)} : \text{the transition from the state} \ q_1 \ \text{to the state} \ q_{(m+1)+1}
                           \  \text{inputing} \  z_{q_1},\\
&\text{(the third term)} : \text{the transition from the state} \ q_1 \ \text{to the state} \ q_{(m+1)+2} 
                           \ \text{inputing} \  z_{p_1+q_1}.
\end{align*}
Thus we have \eqref{fig:h}. \hfill$\square$
\end{pf}
 
\begin{cor}
The harmonic automaton of $w_1^{\star} * w_2^{\star}$ is represented by
the following transtion diagram \eqref{fig:hstar}:

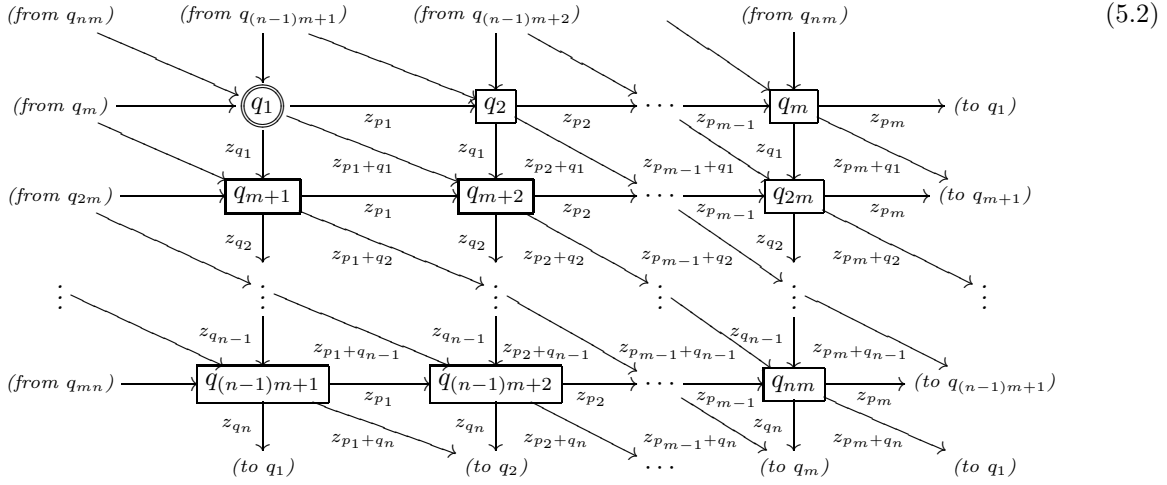
\begin{align}\label{fig:hstar}
\UseTips
\xymatrix @-2mm {
  \text{\scriptsize{(from $q_{nm}$)}} \ar[rd]&
  \text{\scriptsize{(from $q_{(n-1)m+1}$)}} \ar[d] \ar[rd]
& \text{\scriptsize{(from $q_{(n-1)m+2}$)}}\ar[d] \ar[rd]
& \ar[rd]
& \text{\scriptsize{(\text{from $q_{nm}$})}} \ar[d]
&
\\
\text{\scriptsize{(from $q_m$)}} \ar[rd] \ar[r]
& *++[o][F=]{q_1} \ar[r]_{z_{p_1}} \ar[d]_{z_{q_1}}@\ar[rd]_{z_{p_1+q_1}}
& *+[F]{q_2} \ar[r]_{z_{p_2}} \ar[d]_{z_{q_1}} \ar[rd]_{z_{p_2+q_1}}
& \cdots \ar[r]_{z_{p_{m-1}}} \ar[rd]_{z_{p_{m-1}+q_1}}
& *+[F]{q_m} \ar[d]_{z_{q_1}} \ar[r]_{z_{p_m}} \ar[rd]_{z_{p_m+q_1}}
&\text{\scriptsize{(to $q_1$)}}
\\
\text{\scriptsize{(from $q_{2m}$)}} \ar[rd] \ar[r]
& *+[F]{q_{m+1}} \ar[r]_{z_{p_1}} \ar[d]_{z_{q_2}} \ar[rd]_{z_{p_1+q_2}}
& *+[F]{q_{m+2}} \ar[r]_{z_{p_2}} \ar[d]_{z_{q_2}} \ar[rd]_{z_{p_2+q_2}}
& \cdots \ar[r]_{z_{p_{m-1}}} \ar[rd]_{z_{p_{m-1}+q_2}}
& *+[F]{q_{2m}} \ar[d]_{z_{q_2}} \ar[r]_{z_{p_m}} \ar[rd]_{z_{p_m+q_2}}
&\text{\scriptsize{(to $q_{m+1}$)}}
\\
  \vdots  \ar[rd]
& \vdots  \ar[d]_{z_{q_{n-1}}} \ar[rd]_{z_{p_1+q_{n-1}}}
& \vdots  \ar[d]_{z_{q_{n-1}}} \ar[rd]_{z_{p_2+q_{n-1}}}
& \vdots  \ar[rd]_{z_{p_{m-1}+q_{n-1}}}
& \vdots  \ar[d]_{z_{q_{n-1}}} \ar[rd]_{z_{p_m+q_{n-1}}}
& \vdots  
\\
\text{\scriptsize{(from $q_{mn}$)}} \ar[r]
& *+[F]{q_{(n-1)m+1}} \ar[r]_{z_{p_1}} \ar[d]_{z_{q_n}} \ar[rd]_{z_{p_1+q_n}}
& *+[F]{q_{(n-1)m+2}} \ar[r]_{\quad z_{p_2}} \ar[d]_{z_{q_n}} \ar[rd]_{z_{p_2+q_n}}
& \cdots \ar[r]_{z_{p_{m-1}}} \ar[rd]_{z_{p_{m-1}+q_n}}
& *+[F]{q_{nm}} \ar[d]_{z_{q_n}} \ar[r]_{z_{p_m}\quad } \ar[rd]_{z_{p_m+q_n}}
&\text{\scriptsize{(to $q_{(n-1)m+1}$)}}
\\
&
 \text{\scriptsize{(to $q_1$)}}
&\text{\scriptsize{(to $q_2$)}}
&\cdots
&\text{\scriptsize{(to $q_m$)}}
&\text{\scriptsize{(to $q_1$)}}
}
\end{align}
\end{cor}

\newtheorem{hkkpf}{Proof of Theorem \ref{thm:hkk}}
\renewcommand{\thehkkpf}{}

\begin{hkkpf}
We prove the case of $m=3$:
Let $\omega$ be a primitive cubic root of unity. The harmonic automaton representing
$z_k^{\star}*(\omega z_k)^{\star}$ is as follows:

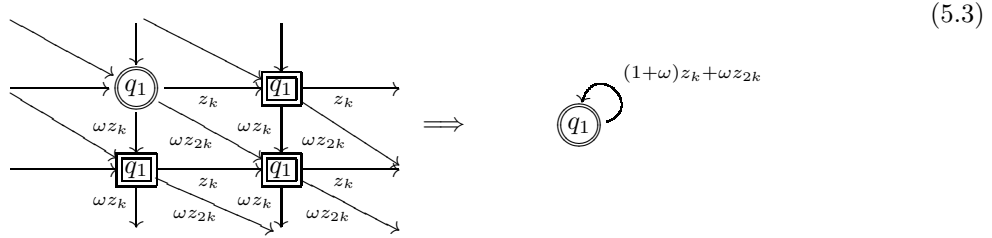
\begin{align}\label{fig:hkkk1}
\UseTips
\xymatrix @-3mm {
   \ar[rrd]
&
&  \ar[d] \ar[rrd]
&
&  \ar[d] 
& &
\\
\ar[rrd] \ar[rr]
&
& *++[o][F=]{q_1} \ar[rr]_{z_k} \ar[d]_{\omega z_k} \ar[rrd]_{\omega z_{2k}}
&
& *+[F=]{q_1} \ar[rr]_{z_k} \ar[d]_{\omega z_k} \ar[rrd]_{\omega z_{2k}}
& &
\\
 \ar[rr]
 &
& *+[F=]{q_1} \ar[rr]_{z_k} \ar[d]_{\omega z_k} \ar[rrd]_{\omega z_{2k}}
&
& *+[F=]{q_1} \ar[rr]_{z_k} \ar[d]_{\omega z_k} \ar[rrd]_{\omega z_{2k}}
& & 
\\
&
&
&
&
& &
}
\UseTips
\xymatrix @+2mm {
 & \\
\Longrightarrow  
& *++[o][F=]{q_1} 
\ar@(r,u)[]_{(1+\omega)z_k+ \omega z_{2k}}
}
\end{align}
Hence we have \
$z_k^{\star}*(\omega z_k)^{\star}=
\{(1+\omega)z_k+ \omega z_{2k}\}^{\star}$. 
The harmonic product of $(z_k^{\star}*(\omega z_k)^{\star})*\omega^2 z_k$ is represented by
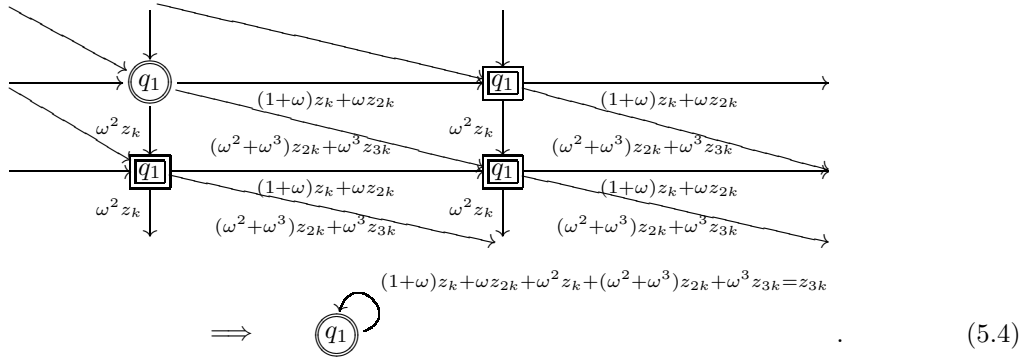
\begin{align}\label{fig:hkkk2}
\UseTips
\xymatrix @-2mm {
   \ar[rrd]
&
&  \ar[d] \ar[rrrrrd]
&
&
&
&
&  \ar[d] 
&
&
&
& &
\\
\ar[rrd] \ar[rr]
&
& *++[o][F=]{q_1} \ar[rrrrr]_{(1+\omega)z_k+ \omega z_{2k}}
 \ar[d]_{\omega^2 z_k}
  \ar[rrrrrd]_{(\omega^2+\omega^3)z_{2k}+ \omega^3 z_{3k}}
  &
  &
&
&
& *+[F=]{q_1} \ar[rrrrr]_{(1+\omega)z_k+ \omega z_{2k}}
 \ar[d]_{\omega^2 z_k}
  \ar[rrrrrd]_{(\omega^2+\omega^3)z_{2k}+ \omega^3 z_{3k}}
  &
  &
&
& &
\\
 \ar[rr]
 &
& *+[F=]{q_1} \ar[rrrrr]_{(1+\omega)z_k+ \omega z_{2k}}
 \ar[d]_{\omega^2 z_k}
  \ar[rrrrrd]_{(\omega^2+\omega^3)z_{2k}+ \omega^3 z_{3k}}
  &
  &
&
&
& *+[F=]{q_1} \ar[rrrrr]_{(1+\omega)z_k+ \omega z_{2k}}
 \ar[d]_{\omega^2 z_k}
  \ar[rrrrrd]_{(\omega^2+\omega^3)z_{2k}+ \omega^3 z_{3k}}
  &
  & 
  & 
& & 
\\
&
&
&
&
&
&
&
&
&
&
& &
} \notag \\ 
\hspace{1cm}
\UseTips
\xymatrix @-2mm {
\Longrightarrow  
& *++[o][F=]{q_1} 
\ar@(r,u)[]_{(1+\omega)z_k+ \omega z_{2k}+\omega^2 z_k+
(\omega^2+\omega^3)z_{2k}+ \omega^3 z_{3k}= z_{3k}}
}.
\end{align}
Therefore the element accepeted by this automaton is $(z_{3k})^{\star}$. \hfill$\square$
\end{hkkpf}

\appendix
\section{On the value of $\zeta(\{2k\}_n)$}

Let $\omega=\omega_{2k}$ be a primitive $2n$-th root of unity. From
\begin{align*}
    \frac{\sin\pi x}{\pi x} = \prod_{r=1}^{\infty} \left( 1 - \frac{x^2}{r^2} \right),
\end{align*}
we have
\begin{align}
  \frac{\sin\pi x \sin\pi\omega x \cdots\cdots \sin\pi\omega^{k-1} x}{\pi^k\omega^{k(k-1)/2}x^k}
  = \prod_{r=1}^{\infty} \left(1 - \frac{x^{2k}}{r^{2k}} \right).
\end{align}
The Taylor expansion of the RHS above is
\begin{align}\label{eq:RHS}
        \sum_{n=0}^{\infty}\zeta(\{2k\}_n)(-x^{2k})^n ,
\end{align}
so it is easy to see 
\begin{align}
       \zeta(\{2\}_n) = \frac{\pi^{2n}}{(2n+1)!}, \quad
       \zeta(\{4\}_n) = \frac{2^{2n+1}\pi^{4n}}{(4n+2)!}.
\end{align}
Comparing the Taylor expansion of the LHS 
\begin{align*}
   & \frac{1}{(2i)^k\pi^k\omega^{k(k-1)/2}x^k}
                   \sum_{\varepsilon_0=\pm1,\ldots,\varepsilon_{k-1}=\pm1}
                      \varepsilon_0\cdots\varepsilon_{k-1}
                        e^{\pi ix(\varepsilon_0+\varepsilon_1\omega
                         +\cdots+\varepsilon_{k-1}\omega^{k-1})} \\
   = & \frac{1}{2^k\omega^{k(k-1)/2}}\sum_{n=0}^{\infty}\frac{(\pi i )^{2nk}}{(2nk+n)!}
                \sum_{\varepsilon_0=\pm1, \ldots , \varepsilon_{k-1}=\pm1}
                 \varepsilon_0\cdots\varepsilon_{k-1}(\varepsilon_0+\varepsilon_1\omega
                            +\cdots+\varepsilon_{k-1}\omega^{k-1})^n
\end{align*}
with \eqref{eq:RHS} we obtain the following proposition:
\begin{prop}We have
\begin{align}
     \zeta(\{2k\}_n) = \frac{(-1)^{(k+1)n}\pi^{2kn}}{2^k\omega^{k(k-1)/2}(2kn+k)!}
     \, \sum_{\varepsilon_0=\pm1, \ldots , \varepsilon_{k-1}=\pm1}
        \varepsilon_0\cdots\varepsilon_{k-1}(\varepsilon_0+\varepsilon_1\omega
                            +\cdots+\varepsilon_{k-1}\omega^{k-1})^{2kn+k}.
\end{align}
\end{prop}
From this proposition, as the explicit formula for $\zeta(\{2k\}_n) \ (k=3,4,5,6)$, we have
\begin{align}
     & \zeta(\{6\}_n) = \frac{6(2\pi)^{6n}}{(6n+3)!} \ ,  \\
     & \zeta(\{8\}_n) = \frac{2^{6n+2}\pi^{8n} \{(3+2\sqrt{2})^{2n+1}+(3-2\sqrt{2})^{2n+1}\}}
                            {(8n+4)!} \ , \\
     & \zeta(\{10\}_n) = \frac{5\cdot2^{8n}\pi^{10n} \{ 2^{2n+1}+
                                     (11+5\sqrt{5})^{2n+1}+(11-5\sqrt{5})^{2n+1} \}}
                            {(10n+5)!} \ , \\
     & \zeta(\{12\}_n) = \frac{3\cdot2^{12n+2}\pi^{12n}\{ 2^{6n+3}+
                     (26+15\sqrt{3})^{2n+1}+(26-15\sqrt{3})^{2n+1} \}}{(12n+6)!} \ .
\end{align}
These formulas reflect the fact that $\mathbb{Q}(\sqrt{2}) \subset \mathbb{Q}(\omega_8), \
\mathbb{Q}(\sqrt{5}) \subset \mathbb{Q}(\omega_{10}), \ 
\mathbb{Q}(\sqrt{3}) \subset \mathbb{Q}(\omega_{12})$. 
To represent $\zeta(\{14\}_n)$, we need the imaginary quadratic field $\mathbb{Q}(\sqrt{-7})$
and the cubic field $\mathbb{Q}(2\cos\frac{2\pi}{7})$ which are subfields of 
$\mathbb{Q}(\omega_{14})$: Let $\alpha_1, \ \beta_1$ be the roots of the quadratic equation
\begin{align*}
     \lambda^2 -13\lambda+128=0,
\end{align*}
and $\alpha_2, \ \beta_2, \ \gamma_2$ be the roots of the cubic equation
\begin{align*}
        \lambda^3-57\lambda^2+103\lambda-1=0.
\end{align*}
Then 
\begin{align}
     \zeta(\{14\}_n)= & \frac{7\cdot2^{14n+1}\pi^{14n}}{(14n+7)!}
       \left\{ 1+\alpha_1^{2n+1}+\beta_1^{2n+1}+\alpha_2^{2n+1}+\beta_2^{2n+1}
                  +\gamma_2^{2n+1}  \right. \nonumber \\
               & \hspace{30mm} \left.  + \alpha_2^{-2n-1}+\beta_2^{-2n-1}+\gamma_2^{-2n-1} \right\}.
\end{align}


\begin{thebibliography}{Mittelbach100}

\bibitem[AK]{AK}
T. Arakawa, M. Kaneko,
\newblock Notes on Multiple Zeta Values and Multiple $L$ Values,
\newblock {\em lecture note delivered at Rikkyou University (in Japanese)},
\newblock (2002).

\bibitem[An]{An}
A. Kelarev,
\newblock GRAPH ALGEBRAS AND AUTOMATA,
{\em Marcel Dekker,Inc. New York $\cdot$ Basel,} (2003).

\bibitem[$\textrm{B}^2$]{BB}
D. Bowman, D. M. Bradley,
\newblock The Algebra and Combinatorics of Shuffles 
and Multiple Zeta Values,
\newblock {\em J. Comb. Theory,
 Series A,} $\textbf{97}$ (2002), 43-61.

\bibitem[$\textrm{B}^3$L]{BBBL}
Jonathan M. Borwein, D. M. Bradley, 
D. J. Broadhurst, P. Lison$\check{\textrm{e}}$k,
\newblock Combinatorial Aspects
 of Multiple Zeta Values,
\newblock {\em Elec. J. Comb.,} 
$\textbf{5}$ (1998), No.1, $\sharp$ R38.

%\bibitem[$\textrm{B}^3$L2]{BBBL2}
%Jonathan M. Borwein, David M. Bradley, 
%David J. Broadhurst, Petr Lison$\check{\textrm{e}}$k,
%\newblock Special Values of Multidimensional Polylogarithms,
%\newblock {\em Trans. Amer. Math.Soc.,} 
%$\textbf{355}$ (2001), No.3, 907-941

\bibitem[H]{H}
M. Hoffman,
\newblock The algebra of multiple harmonic series
\newblock {\em J. of Alg.,} 
$\textbf{194}$ (1997), 477-495.

\bibitem[HMU]{HMU}
J. E. Hopcroft, R. Motwani, J. D. Ullman,
\newblock Introduction to Automata Theory, Languages, and Computation, 
Second Edition,
{\em ADDITION WESLEY LONGMAN, a Pearson Education Company,}
(2001).

\bibitem[K]{K}
S. Kitani,
\newblock Relations of Multiple Zeta Values and Multiple $L$ Values from
the viewpoint of Automaton Theory,
\newblock {\em master thesis at the graduate school of Waseda university (in Japanese)},
(2004).

\bibitem[S]{S}
E. Sawada,
\newblock Automaton and Multiple Zeta Values,
\newblock {\em master thesis at the graduate school of Waseda university (in Japanese)},
(2004).

\bibitem[W]{W}
M. Waldschmidt,
\newblock Multiple Polylogarithms,
\newblock (2001). \\
http://www.math.jussieu.fr/`miw//articles/ps/mpl.ps

\bibitem[Z]{Z}
D. Zagier, 
\newblock Values of zeta functions and their applications, 
\newblock {\em First European Congress of Mathematics 
(Paris,1992), vol 2, Progress in Math,}
$\textbf{120}$ (1994), 497-512.

\end{thebibliography}
\end{document}